\documentstyle[12pt,amssymb,amscd,oldgerm]{amsart}  % amslatex

\newtheorem{thm}{Theorem}[subsection]
\newtheorem{mthm}{Theorem}

\newtheorem{lemma}[thm]{Lemma}

\newtheorem{definition}[thm]{Definition}
\newtheorem{proposition}[thm]{Proposition}

\theoremstyle{remark}

\theoremstyle{definition}

\numberwithin{equation}{section}

\newcommand{\C}{{\Bbb C}}
\newcommand{\bbA}{{\Bbb A}}
\newcommand{\bbC}{{\Bbb C}}

\newcommand{\bbL}{{\Bbb L}}
\newcommand{\bbR}{{\Bbb R}}
\newcommand{\bbW}{{\Bbb W}}
\newcommand{\bbZ}{{\Bbb Z}}

\newcommand{\bbO}{{\Bbb O}}

\newcommand{\cW}{{\cal W}}

\newcommand{\fgl}{{\frak {gl}}}
\newcommand{\fg}{{\frak g}}
\newcommand{\fh}{{\frak h}}
\newcommand{\fsu}{{\frak su}}
\newcommand{\fu}{{\frak u}}

\newcommand{\Hom}{\operatorname{Hom}}

\newcommand{\isomoto}{\overset{\sim}{\to}}
\newcommand{\id}{\operatorname{id}}

\newcommand{\fd}{{\frak{d}}}
\fontfamily{yfrak}
\selectfont

\begin{document}

\title{Deformation quantization of Endomorphism bundles }

\author{ Johannes Aastrup }
\address{Mathematisches Institut Der Universit\"{a}t M\"{u}nster, Einsteinstr. 62, D-48149 M\"{u}nster}
\email{aastrup@@math.uni-muenster.de}

\maketitle

\tableofcontents

\section{Introduction}
A deformation quantization of a Poisson manifold $(M, \{ \cdot, \cdot
\})$ is an associative product $*$ on $C^\infty (M)[[\hbar]]$, so that
$$f*g=fg+ \hbar P_1(f,g)+\hbar^2 P_2 (f,g)+ \ldots ,$$ $*$ being $\hbar$
linear, $P_i$ bidifferential and
$$[f,g]=i\hbar \{ f,g \}+ O(\hbar),$$
where $[\cdot , \cdot ]$ is the commutator with respect to $*$. 

This notion was first introduced in \cite{BFFLS}.

The question of existence and classification of deformation
quantizations on general Poisson manifolds was solved in 1997 by
Kontsevich in \cite{ko}.

The simpler case of existence of deformation quantizations of the
canonical Poisson structure on   symplectic manifolds was solved
already in \cite{WL}. A simple geometric construction of deformation
quantizations of symplectic manifolds were given by Fedosov in
\cite{Fe3}. The advantage of Fedosov's construction compared to the
ones in \cite{ko} and \cite{WL} is that it is easy to handle and also
suitably generalizable. The most general  setting of the Fedosov
construction is probably given in \cite{NT4}, where deformation
quantizations of symplectic Lie algebroids is done. Also the
classification of deformation quantizations becomes amenable in view
of the Fedosov construction. In the case of symplectic manifolds this
was done in \cite{NT1} and the classification of deformation
quantizations on a symplectic manifold $(M, \omega )$ is given by the
points (characteristic classes) $\theta$ in the space
$$\frac{\omega}{i \hbar }+H^2(M, \bbC [[\hbar ]]).$$

One of the main class of examples of deformation quantizations of symplectic manifolds
are those coming from asymptotic calculus of pseudodifferential
operators on manifolds, see for example \cite{NT7}. If we consider the
asymptotic
calculus of
pseudodifferential operators on a manifold $M$, we will get a
deformation quantization of the cotangent bundle $T^*M$ of $M$, where
$T^*M$ is equipped with the canonical symplectic structure.

This example gives the connection to index theory. On a deformation
quantization of any symplectic $2n$ dimensional manifold there is a
canonical trace, unique up to multiplication by a scalar,
of the form
\begin{equation} \label{sp}
Tr(a)=\frac{1}{n!(i \hbar)^n}\left( \int_Ma\omega^n + O(\hbar)\right).
\end{equation}
By an appropriate choice of the representation of the quantization,
i.e. after applying a linear isomorphism of $C^\infty (M)[[\hbar]]$ of
the form
$$
f\rightarrow f+\hbar D_1(f)+\ldots
$$
one can assure that $Tr$ has the form
$$
Tr(a)=\frac{1}{n!(i \hbar)^n} \int_Ma\omega^n ,
$$
which fixes it uniquely.

In most proofs of the Atiyah-Singer index theorem and related
``local'' index theorems, one of the main difficulties is to compute
the trace of a certain operator on a Hilbert space, usually
$L^2(M)$ as above. In order to compute this trace a scaling $\hbar$ in $\bbR_+$
of the operator is introduced, and the asymptotic expansion of the
trace as $\hbar \rightarrow 0$ becomes computable, at least the
constant term in the expansion. The computations coming out of this is
computations like \ref{sp}. This is why computing the canonical trace
on deformation quantizations is called algebraic index
theory. Actually,  according to \cite{NT7}, computing the trace on deformation quantizations,
 in a way that will be described now,
implies the Atiyah-Singer index theorem.

Many elements, though not all, on which computing the trace is interesting,
are first components in classes in cyclic periodic homology. The
cyclic periodic homology or rather cohohmology was invented by Connes
in \cite{Co1}. It is the noncommutative analog of De Rahm cohomology
and was already at the beginning intimately connected to index
theory. A complex computing the cyclic periodic homology of a unital
algebra $A$ over a field $k$ is given by
$$
CC_{even}^{per}(A)=\prod_i A\otimes \bar{A}^{\otimes 2i}; \quad
CC_{odd}^{per}(A)=\prod _i A\otimes \bar{A}^{\otimes 2i+1}
$$
where $\bar{A}=A/k\cdot 1$
and the differential
$$
CC_{even}^{per}(A)
\stackrel{b+B}{\longleftrightarrow}CC_{odd}^{per}(A)
$$
is given by
$$
\begin{array}{l}
b(a_0 \otimes \ldots \otimes a_n)= \\
\sum_{k=0}^{n-1}(-1)^ka_0\otimes
  \ldots \otimes
  a_k a_{k+1} \otimes \ldots \otimes a_n
+ (-1)^{n}a_na_0 \otimes a_1
  \otimes \ldots \otimes a_{n-1}
\end{array}
$$
and
$$
B(a_0 \otimes \ldots \otimes a_n)=\sum_{k=0}^n (-1)^k 1\otimes
a_k\otimes a_{k+1} \ldots \otimes a_n \otimes a_0 \ldots a_{k-1}.
$$
If for example $p \in M_n(A)$ is a projection, a class in
$HC^{per}_{even}(A)$, the Chern character of p is given by the formula
$$
tr (p+ \sum_{k \geq 1}\frac{(2k)!}{k!} (p-1/2)\otimes p^{\otimes 2k}),
$$
where
$$
tr:M_n(A)^{\otimes k}\rightarrow A^{\otimes k}
$$
is the
map given by
$$
(M_1\otimes a_1) \otimes \ldots \otimes (M_k\otimes
a_k) \mapsto tr(M_1\ldots M_k)a_1\otimes \ldots \otimes a_k
$$
Therefore $tr (p)$ can be regarded as the first component of a class in
cyclic periodic homology.

Evaluating $Tr$ on the first component gives a morphism of complexes
\begin{equation}\label{spfunk}
Tr: CC_*^{per} (A^\hbar_c) \rightarrow \bbC [[\hbar , \hbar^{-1}],
\end{equation}
where $A^\hbar_c$ is the algebra of compactly supported elements
$C_c^\infty (M)[[\hbar]])$ in a deformation
quantization $A^\hbar$, and $ \bbC [[\hbar , \hbar^{-1}]$ is
considered as a complex concentrated at degree zero with the trivial
boundary map.

Computing the trace on elements that are first components of a class in
cyclic periodic homology is therefore the same as computing \ref{spfunk} at the
level of homology.

In \cite{NT1} it is proved that
$$ Tr ( \cdot ) \sim (-1)^n \int_M \hat{A}(TM)e^{\theta}\tilde{\mu}(\cdot
)$$
where $\sim$ means that the two sides define the same  morphism at the level
of homology. Here $\theta$ is the characteristic class of the deformation and
$\tilde{\mu}$ is the map $CC^{per}_* (A^\hbar ) \rightarrow
\Omega^*(M)$ given by
$$\tilde{\mu} (a_0 \otimes  \ldots \otimes a_k )=  \frac{1}{k!}\tilde{a}_0d
\tilde{a}_1 \cdots d \tilde{a}_k, \quad \tilde{a}_i=a_i \hbox{ mod } \hbar$$

This settles the problem of computing the trace at the level of
homology for deformation quantizations of symplectic manifolds.

\subsection{Contents of the Paper} Below I propose a definition
of deformation quantization of endomorphism bundles over a symplectic
manifold. The motivation is clear: Deformation quantizations of the
trivial line bundle is the algebraic analog of pseudodifferential
operators in line bundles and therefore deformation quantization of
an endomorphism bundle $End(E)$, $E$ vector bundle over $M$, should be
the algebraic analog of pseudodifferential operators in any vector
bundle having $End(E)$ as endomorphism bundle.

The definition proposed  requires a product $*$
on $$\Gamma (End(E))[[\hbar ]],$$ so the algebra $(\Gamma (End(E))[[\hbar
]], *)$ is locally isomorphic to $M_N(\cW_n)$, i.e. $\cW_n$ being the Weyl
algebra, the canonical deformation quantization of the standard
symplectic structure on $\bbR^{2n}$.

It turns out that the Fedosov construction also works in this case.
Thus we can  let $M_N(\bbA^\hbar )$ be the algebra of jets at zero of
elements in $M_N(\cW_n)$.  Associated to $(M, \omega, End(E))$ there is
an algebra bundle $\bbW$ with fiber $M_N(\bbA^\hbar )$. Put $\fg
=Der(M_N(\bbA^\hbar ))$. There is now a short exact sequence
$$0 \rightarrow \frac{1}{\hbar} \bbC [[\hbar ]] \rightarrow \tilde{\fg
  }\rightarrow \fg \rightarrow 0$$
of Lie algebras.
 The Fedosov construction then consist, for a given element
$\theta$  in $$ \frac{\omega}{i \hbar }+H^2(M,\bbC[[\hbar ]]),$$
in constructing a flat connection $\nabla $ in $\bbW$ with values in
$\fg$, such that $ker \nabla \simeq \Gamma (End(E))[[\hbar ]]$ linearly
and $\nabla$ admits a lift $\tilde{\nabla}$ to a connection with
values in $\tilde{\fg}$ and curvature $\theta$. The product on $\Gamma
(\bbW)$ induces a product on $ker \nabla =\Gamma (End(E))[[\hbar ]]$.
This product gives a deformation quantization of $End(E)$ and $\theta$
will be an isomorphism invariant of the deformation quantization.

This construction is done in section \ref{defqe}. In this
section the following is also shown.
\begin{mthm}
A deformation quantization of
$End(E)$ is isomorphic to the flat sections of a Fedosov
connection, and the isomorphism classes of deformation quantizations
of
$End(E)$ are
classified by the points in
$$
\frac{\omega}{i \hbar }+H^2(M, \bbC [[\hbar ]]).
$$
\end{mthm}
The principle,  that a deformation quantization comes as flat sections
in a certain infinite dimensional vector bundle, is not special to
deformation quantizations. In section \ref{end} it is shown, that
sections of $\Gamma (End(E))$ are flat sections in an algebra bundle
with fibre $$M_N(\bbC [[\hat{x}_1, \ldots , \hat{x}_n ]]).$$

The reason
for redoing this construction for $End(E)$ is, that it is notationally
simpler, and, hopefully,
clarifies the construction. Therefore in section \ref{defqe} only the
differences in the construction for $End(E)$ and deformation
quantizations of endomorpism bundles are spelled out.

Like in the scalar case there are canonical traces on deformation
quantizations of endomorphism bundles. The rest of the paper is
devoted to an index theory for these traces.  The methods used for this
have been developed by Nest and Tsygan in \cite{NT1},
\cite{NT5}, \cite{BNT} and \cite{NT4}. These methods are based on the following.

\begin{enumerate}
\item The action of the reduced cyclic
  complex  $\overline{C}^\lambda_*(A)$, on $CC^{per}_*(A)$:
$$
 \overline{C}^\lambda_*(A) \times CC^{per}_*(A)
 \rightarrow CC^{per}_*(A)
$$
\item The construction of the fundamental class, a special class  in 
$\overline{C}^\lambda_*(A^\hbar_E)$, where
  $A^\hbar_E$ is the deformation quantization of $End(E)$. Or
  rather the construction of a class in the \v{C}ech complex,  $\check{C}^*(M,
  \overline{C}^\lambda_*(A^\hbar_E))$, with
  values in the presheaf $V\rightarrow
  \overline{C}^\lambda_*(A^\hbar_{E|V})$.

\item Computations in Lie algebra cohomology in order to identify the
  fundamental class at the level of cohomology.
\end{enumerate}

\noindent The fundamental class $U$  lives in
$\overline{C}^\lambda_{2n-1}(M_N(\cW_n))$. Its role is that it
relates $Tr$ to $\tilde{\mu }$ when evaluated at
  classes that are scalar mod $\hbar$. This has the effect that
$$
Tr(U \cdot a)=(-1)^n\int \tilde{\mu}(a), \quad a \in
CC^{per}_*(M_N(\cW_{n,c})).
$$
Here, as before, the subscript ``c'' denotes the ideal of compactly
supported elements in the deformation algebra in question.

The following plays the major role.
\begin{mthm}
In cohomology, the class $U$ has a unique
extension to a class in
$\check{C}^*(M,\overline{C}^\lambda_*(A^\hbar_E))$, also denoted by $U$.
On classes of the form  $a_0 \otimes \ldots \otimes a_k \in
CC^{per}_*(A^\hbar_{E,c})$, where $a_i$ is a scalar mod $\hbar$, the following holds:
$$
Tr(U \cdot a)=(-1)^n \int \tilde{\mu}(a) \hbox{ mod }\hbar \quad .
$$
\end{mthm}
It is not difficult to see that this implies for general classes $a_0
\otimes \ldots \otimes a_k$ in
$CC^{per}_*(A^\hbar_{E,c})$ that
$$Tr(U \cdot a_0\otimes \ldots \otimes a_k)= (-1)^n\int
ch(End(E))^{-1}ch(\nabla )(\tilde{a}_0\otimes \ldots \otimes
\tilde{a}_k) \hbox{ mod }\hbar$$
where $ch(End(E))$ is the usual Chern character of $End(E)$ as a
vector bundle, $\nabla$ is a connection in $End(E)$ and
\begin{eqnarray*} && ch(\nabla)(\tilde{a}_0 \otimes \ldots \otimes
  \tilde{a}_k)= \\
&& \int_{\Delta_k}tr (\tilde{a}_0e^{-t_0\nabla^2}\nabla
(\tilde{a}_i)e^{-t_1\nabla^2}\cdots \nabla
(\tilde{a}_k)e^{-t_k\nabla^2} )dt_0 \cdots dt_{k-1}
\end{eqnarray*}
the J.L.O. cocycle  associated to $\nabla$.

To finish the index theory, the fundamental class has to be
identified. This is done via Lie algebra cohomology. We have the Gelfand-Fuks
morphism of complexes
$$C^*(\fg, \fsu(N)+\fu(n); \overline{C}^\lambda_*(M_N(\bbA^\hbar ))) \rightarrow
\check{C}^*(M, \Omega^*(M, \overline{C}^\lambda_*(M_N(\bbA^\hbar)))),$$
the latter complex being quasi isomorphic to
$\check{C}^*(M,\overline{C}^\lambda_*(A^\hbar_E))$.  As in the case of
$M_N(\cW_n)$, there is a fundamental class $U$ in
$\overline{C}^\lambda_*(M_N(\bbA^\hbar))$ extending uniquely in
cohomology to a class in Lie algebra cohomology, also denoted by $U$. It
turns out that $GF(U)$ is equivalent to $U$ in
$\check{C}^*(M,\overline{C}^\lambda_*(A^\hbar_E))$ via the quasi-isomorphism between $\check{C}^*(M,\overline{C}^\lambda_*(A^\hbar_E))$
and $\check{C}^*(M, \Omega^*(M,
\overline{C}^\lambda_*(M_N(\bbA^\hbar))))$. Hence the question of
computing or identifying $U$ is now a question of computations in
$$C^*(\fg, \fsu(N)+\fu(n); \overline{C}^\lambda_*(M_N(\bbA^\hbar
))).$$
It turns out to be useful to work with the differential graded
algebra $M_N(\bbA^\hbar)[\eta ]$, where $\eta$ is a formal variable,
$\eta^2=0$ and the differential is given by $\frac{\partial}{\partial
  \eta }$. 

The reason for doing this is to include the
identity operation on $CC^{per}_*(M_N(\bbA^\hbar))$. Thus the action
of $\overline{C}^\lambda_*(M_N(\bbA^\hbar))$ on
$CC^{per}_*(M_N(\bbA^\hbar))$ extends to an action of
$\overline{C}^\lambda_*(M_N(\bbA^\hbar)[\eta ])$  on
$CC^{per}_*(M_N(\bbA^\hbar))$. With this action the classes
$\eta^{(k+1)}=k!\eta^{\otimes k+1}$ becomes the identity
operations. The main technical theorem of this paper, Theorem
\ref{teks}, states the following
\begin{mthm}
$U$ is equivalent to
$$\sum_{m \geq 0}(\hat{A}\cdot e^\theta \cdot (ch)^{-1})^{-1}_{2m}\cdot
\eta^{(m)}$$
in $C^*(\fg, \fsu(N)+\fu(n); \overline{C}^\lambda_*(M_N(\bbA^\hbar
)[\eta]))$,
i.e. it defines the same cohomology class. 

Here $\hat{A}$ is the Lie algebraic
$\hat{A}$ coming from $\fu (n)$, $ch$ is the Lie algebraic Chern
character coming from $\fsu (N)$ and $\theta$ is the Lie algebraic
class of deformation. 
\end{mthm}
From this follows the index theorem

\begin{mthm}
Let $A^\hbar_E$ be a deformation quantization of $End(E)$ and let $\theta$ be the associated characteristic class. Then the identity 
$$Tr(a)=(-1)^n\int \hat{A}\cdot e^\theta \cdot ch(\nabla)(\tilde{a})$$
holds when $a$ is a cycle in $CC^{per, \bbC}_*(A^\hbar_{E,c})$.  
\end{mthm}

\textit{Acknowledgements.} I would like, first of all, to thank my advisor Ryszard Nest. I would also like to thank  Boris Tsygan for
helpful conversations and Paulo Almeida and Nuno Martins for having
hosted me for two months at the I.S.T. in Lisbon.

%\section{Structure of $Aut(M_N(C^\infty(\bbR^n)))$ and $Der(M_N(C^\infty(\bbR^n)))$}

%\begin{lemma} Derivations of $M_N(C^\infty(\bbR^n))$ are on the form
%$X+[a,\cdot]$, where $X$ is a vector field and $a$ is in
%$M_N(C^\infty(\bbR^n))$.
%\end{lemma}

%\textit{Proof.} Consider the family of derivations of $M_N(\bbC)$ given by
%$$x\rightarrow D(A)(x)$$
%Since the derivations are inner we get an element $a\in
%M_N(C^\infty(\bbR^n))$ such that $D(A)(x)=[a(x),A]$. We thus see,
%that for $E=D-[a,\cdot]$ we have
%$E(a)$ is scalar valued when $a$ is scalar valued and hence given by
%a vector field. We thus get, the desired.
%\\

%We now look at the structure of $Aut(M_N(C^\infty(\bbR)))$. Given an
%automorphism $\Phi$ of $M_N(C^\infty)$,  this will also be an
%automorphism of the center and hence given  by a
%diffeomorphism. Applying the inverse of this diffeomorphism on $\phi$
%we get a $C^\infty(\bbR^n)$ module automorphism $\Xi$ of
%$M_N(C^\infty(\bbR))$. But  then $\xi$ must be given by a smooth
%family of automorphisms of $M_N(\bbC)$. But an automorphism of this
%form is given by conjugating with an invertible element in
%$M_N(C^\infty(\bbR^n))$.  We thus have the following
%\begin{lemma} An automorphism of $M_N(C^\infty(\bbR^n))$ is given
%conjugating with a diffeomorphism composed with conjugating with an
%invertible element in  $M_N(C^\infty(\bbR^n))$
%\end{lemma}

\section{Sections of Endomorphism Bundles as flat Sections in a profinite Bundle}

\label{end}
Let $\bbO_n=\bbC[[\hat{x}_1,\ldots,\hat{x}_n]]$ and $M_N(\bbO_n )=M_N
\otimes \bbO_n$, where $M_N$ is the algebra of $N \times N$ complex matrices. We
give  $M_N(\bbO_n )$ a grading by setting
$$ \hbox{deg} (b \otimes \hat{x}_i )=1, \quad b \in M_N.$$
Furthermore we give $M_N(\bbO_n)$ the $I$-adic topology, where $I$ is
the ideal generated by elements of degree $\geq 1$.

\begin{definition} Let $G$ be the group of continuous automorphisms of
  $M_N(\bbO_n)$ such that
 the induced automorphism on the
centre $\bbO_n$ is an automorphism induced by  an automorphism on 
$\bbR[[\hat{x}_1,\ldots,\hat{x}_n]]$.
\end{definition}

\begin{lemma}\label{finner} An automorphism in $G$ is the composition of an automorphism induced by an automorphism of $\bbR[[\hat{x}_1,\ldots,\hat{x}_n]]$ and an inner automorphism.
\end{lemma}

\textit{Proof.} Given an automorphism $\Phi$ in $G$ let $\varphi$ be the induced automorphism on $\bbO_n$. Considering $\chi =\Phi\circ (\varphi \otimes id)^{-1}$ we have that $\chi$ is an $\bbO_n$-module map.

For $A$ in
 $M_N(\bbC)$ we have

$$\chi(A)=D_0(A)+\hbox{higher order terms , }D_0(A)\in M_N.$$
Since $D_0$ is an automorphism of $M_N$, it is inner and hence extends
to $M_N(\bbO_n)$. Let $\chi_1=\chi \circ D_0^{-1}$. Then
$$\chi_1(A)=A+\sum_i\hat{x}_iD_i(A)+\hbox{ higher order terms }$$
Each $D_i$ is a derivation of $M_N$ and hence given by a commutator by  an element $B_i$ in $M_N$. Therefore
$$\chi_1 \circ \exp(\sum_i ad(\hat{x}_iB_i))(A)=A+\hbox{ terms of order }\geq 2.$$
Continuing by induction we get a sequence of elements $C_k$ in
$M_N(\bbO_n)$, where $\hbox{deg}(C_k)=k$, with
$$(\chi \circ \exp (ad(C_0)\circ \ldots \circ \exp (ad (C_k))(A)=A +
\hbox{terms of order } \geq k+1.$$
Since the product
$$ \exp (ad(C_0))\circ \ldots \circ \exp (ad (C_k))\circ \ldots$$
converges, we see, by the Hausdorff-Campbell formula, that $\chi$ is inner.
\\

From Lemma \ref{finner} we get
\begin{proposition} The Lie algebra of the Lie-group $G$ is
$$\bbW_n^0 \ltimes (\fgl_N(\bbO_n) / centre ),$$  where $\bbW_n^0$ is the Lie-algebra of formal vector fields vanishing at zero.
\end{proposition}

We note that $Der (M_N(\bbO_n))$ is larger than
$\bbW_n^0 \ltimes (\fgl (\bbO_n) / centre )$, namely
$$\bbW_n \ltimes (\fgl_N(\bbO_n)/centre )$$ where $\bbW_n$ are all formal vector fields  on $\bbR^n$.

\begin{lemma} \label{cont}
Let $X$ be a contractible open subset of $\bbR^k$, let $Aut (M_n)$ be
the automorphism group of $M_n$ and let $Der (M_n)$ be the derivations
of $M_n$. Then any smooth maps
$\varphi_1 :X \rightarrow Aut(M_n)$ and $\varphi_2:X  \rightarrow
Der(M_n)$ lift to smooth maps $\tilde{\varphi}_1 :X \rightarrow Gl_n$
and $\tilde{\varphi}_2: X \rightarrow M_n$.
\end{lemma}

\textit{Proof.} First the case of $\varphi_1$:

Let $\{e_{ij} \}$ be the standard matrix units. The families
$$ x \rightarrow  \varphi_1(x) (e_{ii}) $$ of
projections over $X$ give rise to a family of line bundles $\{ l_i
\}$ over $X$ by
$$ l_i=(\varphi_1(x)(e_{ii}))(\bbR^n).$$
Since $X$ is contractible, these line bundles are trivial. Let $v_1$
be a smooth nowhere vanishing section of $l_1$ and let
$$v_i(x)=(\varphi_1(x)(e_{i1}))(v_1 (x)).$$
Put
$$(\tilde{\varphi}_1 (x))(e_i)=v_i (x)$$
where $e_i$ is the vector  $(a_1, \ldots , a_n) \in \bbR^n$  with
$a_j =0$, $j\not=i$ and $a_i=1$.
Since
$$(\tilde{\varphi}_1 (x)e_{ij}(\tilde{\varphi}_1(x))^{-1})v_k(x)=
\tilde{\varphi}_1(x)e_{ij}e_k=\tilde{\varphi}_1(x) e_i
\delta_{jk}=\delta_{jk}v_i(x)$$
and
\begin{eqnarray*}
& \varphi_1(x)(e_{ij})v_k(x)= (\varphi_1 (x)(e_{ij}))(
\varphi_1(x)(e_{k1}))v_1(x)\\
&= \delta_{kj}
\varphi_1(x)(e_{i1})(v_1(x))=\delta_{kj} v_i(x)
\end{eqnarray*}
we have that $\tilde{\varphi}_1:X \rightarrow Gl_n$ is smooth and that
$$\tilde{\varphi}_1(x)A (\tilde{\varphi}_1(x))^{-1}= \varphi_1(x)(A),
\quad A \in M_n.$$

In the case of $\varphi_2$, let $\varphi : X \times \bbR \rightarrow
Aut(M_n)$ be defined by
$$\varphi (x,t) = \exp (t \varphi_2 (x)).$$
This is clearly a smooth map and hence by the first part of the lemma we get
a smooth lifting of $\varphi$ to $\tilde{\varphi} : X \times \bbR
\rightarrow Gl_n$.

Defining
$$\tilde{\varphi}_2(x)= \frac{\partial}{\partial
  t}(\tilde{\varphi})(x,0)$$
the lemma follows.
\\

In view of lemma \ref{cont} the proof of lemma \ref{finner} actually shows
\begin{lemma} \label{glatloft}
A smooth family of elements in $G$ over a contractible open subset $X$
of $\bbR^k$ lifts to a smooth family over $X$ of elements in the group
of invertible elements in $M_n ( \bbO_n)$.
\end{lemma}

\subsection{Jets of Sections of Endomorphism Bundles}

Let $A$ be an algebra bundle  over a manifold $M$. In this case we let
$$I_m=\{a\in \Gamma (A)|a(m)=0\} .$$
Given this, we define
$$J_mA=\lim_{\stackrel{\leftarrow}{k}}\Gamma (A)/I^k_m.$$
Denote by $J_m$ the quotient map from $\Gamma (A)$ into
$J_mA$. In the following we are only interested in the case where
$A=End(E)$, $E$ a vector bundle. Note that $J_mA\simeq M_N(\bbO_n)$ by
choosing a trivialization; and that any other trivialization leads to
an automorphism of $M_N(\bbO_n)$ belonging to $G$.

If we are given a smooth path of automorphisms $\Phi_t$ in
$G$, we can, according to Lemma \ref{glatloft}, write it
as $\Phi_t=\varphi_t\circ \chi_t$, where $\varphi_t$ is a smooth path
of automorphisms induced by automorphisms of $\bbR[[\hat{x}_1, \ldots
, \hat{x}_n]]$ and $\chi_t$ is a smooth path of inner automorphisms of
$M_N(\bbO_n)$. It is well known that $\varphi_t$ lifts to a smooth
path of local diffeomorphisms $\tilde{\Phi}_t$ of $\bbR^n$ preserving
zero and since $\chi_t$ lifts to a smooth path of invertible elements
in $M_N(\bbO_n)$, it lifts, by the Borel lemma, to a smooth path of invertible elements in $M_N(C^\infty (U))$, where $U$ is an open subset of $\bbR^n$ containing zero. We thus get

\begin{proposition} \label{vectcom} Any smooth path of automorphism in $G$ lifts to a smooth path of local bundle automorphism of $M_N(C^\infty (\bbR^n))$ preserving zero.
\end{proposition}

\subsection{The Frame Bundle}
For a manifold $M$ with a vector bundle $E$ we  define the following
\begin{definition} The frame bundle $\tilde{M}_E$ is given by
$$\tilde{M}_E=\{ (m,\Phi)|m \in M, \Phi : M_N(\bbO_n) \isomoto J_m End(E)\}$$
\end{definition}
We note that $\tilde{M}_E$ is a profinite manifold and in fact a principal bundle over $M$ with fibre $G$.

\begin{proposition} For all $(m, \Phi) \in \tilde{M}_E$ there is an isomorphism
$$\omega :T_{(m,\Phi)}\tilde{M}_E \rightarrow Der(M_N(\bbO_n))$$  satisfying
\begin{eqnarray*}
(i)\quad& \omega (A^*)=A & A \in \bbW^0_n \ltimes ( \fgl _N(\bbO_n)/centre)  \\
(ii)\quad& \varphi^*\omega=ad\varphi^{-1}\omega & \hbox{ for } \varphi \in G \\
(iii)\quad& d\omega + \frac{1}{2}[\omega , \omega ]=0 
\end{eqnarray*}
where $A^*$ is the fundamental vector field corresponding to $A$.

In other words, $\omega$ is a flat connection in $\tilde{M}_E$ with values in $Der(M_N(\bbO_n))$.
\end{proposition}

\textit{Proof.}  Suppose we are given a path in $\tilde{M}_E$, so $\gamma (t)\in \tilde{M}_E$ with $\dot{\gamma}(0)=v$, $\gamma(0)=(m,\Phi)$ and $\gamma (t)=(m_t,\Phi_t)$. This lifts to a path of trivializations $\tilde{\gamma }=(m_t,\tilde{\Phi}_t)$, $\tilde{\Phi}_t:M_N(C^\infty (U)) \rightarrow \Gamma (End(E))$, that maps $0$ to $m_t$.

Define $\omega(v)$ to be the derivation
$$ J_0(a) \rightarrow J_0\left( \frac{d}{dt} (\tilde{\Phi}^{-1} \circ \tilde{\Phi}_t(a)) \right), \qquad a \in M_N(\C^\infty(U))$$
This does not depend on the choice of $\gamma$ or the lifting to $\tilde{\gamma}$ and will be an isomorphism.

 Identity $(i)$ follows, since $\omega$ is the canonical one form on the fibres of $\tilde{M}_E$.

For  identity $(ii)$ we have to compute $\omega (\varphi_*(v))$,
$\varphi\in G$ , but
$$\omega(\varphi_*(v))=\left( J_0(a) \rightarrow J_0\left( \frac{d}{dt}(\tilde{ \varphi}^{-1} \circ \tilde{\Phi}^{-1} \circ \tilde{\Phi}_t \circ \tilde{\varphi} (a))\right) \right) =ad (\varphi^{-1}) \omega (v)$$

 Identity  $(iii)$ is equivalent to
$$\omega ([\omega^{-1}(X),\omega^{-1}(Y) ])=[X,Y], \qquad X,Y \in Der (M_N(\bbO_n))$$
This is actually a consequence of identity $(i)$, because the statement is obvious for $X,Y$ of the form $\frac{\partial}{\partial \hat{x}_i} \in Der(M_N(\bbO_n))$.  For $X,Y$ in $\bbW_n^0 \ltimes (\fgl_N(\bbO_n)/centre)$ it follows because $\omega$ is the canonical one form on the fibres of $\tilde{M}_E$. Hence it suffices to check the case when $X$ is of the form $\frac{\partial}{\partial \hat{x}_i}$ and $Y \in \bbW_n^0 \ltimes (\fgl_N(\bbO_n)/centre)$. Let therefore  $\varphi_t$ be the one parameter group for $Y$ on $G$. We  have
\begin{eqnarray*}
& \omega ([\omega^{-1}(X),\omega^{-1}(Y)])=\omega (\lim_{t \rightarrow 0}\frac{1}{t}((\varphi_t)_* \omega^{-1}(X)-\omega^{-1}(X))) \\
& =\lim_{t \rightarrow 0} \frac{1}{t}(ad(\varphi_t^{-1})(X)-X)=[X,Y],
\end{eqnarray*}
from which the proposition follows.
\\

Given $\tilde{M}_E$, we define the jet bundle of $End(E)$ by
$$JE=\tilde{M}_E\times_G M_N(\bbO_n).$$

The flat connection on $\tilde{M}_E$ gives a flat connection in $JE$
in the following way: Choose a trivialization of $JE_{|U} \isomoto
U\times M_N(\bbO_n)$. This corresponds to a lift $\sigma :U
\rightarrow \tilde{M}_{E|_U}$ of the projection $P:\tilde{M}_E
\rightarrow M$. In this trivialization the connection $\nabla$ is
given by $\nabla =d - \sigma^*(\omega)$, where $\omega$ is the
connection described in proposition \ref{vectcom}.
\begin{proposition} The complex $(\Omega^*(M,JE),\nabla)$ is acyclic and the cohomology is isomorphic to $\Gamma (End(E))$.
\end{proposition}

\textit{Proof.} There is an injective map $j$ from $\Gamma(End(E))$ into $JE$ given by
$$j(\gamma )=((m,\Phi),\Phi^{-1} (J_m\gamma )), \quad \gamma \in
\Gamma (End(E)).$$
 To see that the image of $j$ belongs to the kernel of
$\nabla$  choose a  trivialization in the sense of a
local bundle map $\Phi:M_N(C^\infty (U)) \rightarrow \Gamma (End(E)) $,
$U$ open subset of $\bbR^n$.  Denote the induced map from
$U$ to $M$ by $\Phi$. From this trivialization we get a special
trivialization of $\tilde{M}_E$ by letting $\tilde{\Phi}_u$ denote the
map $M_N(\bbO_n )=J_u(M_N(C^\infty (U))) \rightarrow
J_{\Phi (u)}(\Gamma (End(E)))$ induced from $\Phi$ and then define 
$$ U\times G \ni (u,g)\rightarrow ( \Phi(u), g \tilde{\Phi}_u).  $$

Using this trivialization, we get a local bundle isomorphism
$$C^\infty (U)\otimes M_N(\bbO_n) \rightarrow JE, $$
 and in this trivialization it
is not difficult to see that
\begin{equation} \label{kanon}
\nabla = d - \sum_i dx_i \otimes
  \frac{\partial}{\partial \hat{x}_i}.
\end{equation}
If $\gamma \in M_N(C^\infty (U))$, we have that $j (\gamma)$ is
just given by the Taylor expansion in each point, i.e.
$$ j(\gamma )(u)=\sum_I \frac{\partial^{|I|} \gamma }{\partial
  x^I}\hat{x}^I$$
where $I$ runs through all multi-indices. Hence $j(\gamma)\in \ker (
\nabla )$.

A computation in $(\Omega^*(U, M_N(\bbO_n)), \nabla)$, where $\nabla$
is as in \ref{kanon}, gives that $(\Omega^*(U, M_N(\bbO_n)), \nabla)$
is acyclic and the cohomology is $j(M_N(C^\infty (U))$.

We have thus seen that $(\Omega^*(M,JE),\nabla)$ is locally acyclic
and the cohomology is locally isomorphic to $\Gamma (End(E))$. Since we
 also have seen, that $\Gamma(End(E))$ is globally  contained in $\ker
\nabla$, i.e. the cohomology of $(\Omega^*(M,JE), \nabla)$ is a module over $\Gamma (End(E))$, the statement follows.

\section{Deformation Quantization of Endomorphism Bundles}

\label{defqe}
We start with looking at $\bbR^{2n}$ with the standard symplectic structure and denote the coordinates by
$x_1, \ldots , x_n, \xi_1, \ldots , \xi_n$. On $C^\infty (\bbR^{2n})[[\hbar ]]$ we consider the Weyl quantization given by the product
$$(f*g)(x,\xi)=\exp \left( \frac{i\hbar }{2} \sum_{k=1}^n (\partial_{x_k} \partial _{\eta_k}-\partial_{\xi_k}\partial_{y_k})\right) f(x,\xi )g(y,\eta )|_{(x=y,\xi =\eta )}$$
We will denote the Weyl quantization by $\cW_n$.

Since the definition of the product in the Weyl quantization of two
functions $f,g$ only uses derivatives of $f$ and $g$, the Weyl quantization
makes sense over any open subset $U$ of $\bbR^{2n}$. We will in
this case talk about the Weyl quantization over $U$.

Let  $(M,\omega )$ be a symplectic manifold and let $E$ be a vector bundle over $M$.
\begin{definition}
A deformation quantization of $End(E)$ is a $\hbar$-linear associative product $*$ on $\Gamma (End(E))[[\hbar ]]$, continuous in the $\hbar$-adic topology and satisfying
$$ f*g=fg+\hbar B_i(f,g)+\ldots $$
where $f,g \in \Gamma (End(E))$ and the $B_i$ are bidifferential
expressions. Furthermore we require that $(\Gamma(End(E))[[\hbar ]], *)$ is locally isomorphic to $M_N(\cW_n)$, where $\cW_n$ is the Weyl algebra on some open subset of $\bbR^{2n}$.

In this case locally isomorphic means that we are given a local bundle isomorphism $\Phi: M_N(C^\infty (U))[[\hbar]] \rightarrow \Gamma(End(E))[[\hbar]]|_U$ over a local symplectomorphism, such that the product $*'$, induced by $*$, on $M_N(C^\infty(U))[[\hbar]]$ is isomorphic to $M_N(\cW_n)$, in the sense that there exist differential operators $\{ D_i \}$, such that the map
$$\varphi : (M_N(C^\infty (U))[[\hbar]],*') \rightarrow M_N(\cW_n)$$given by
$$\varphi(a)=a+\hbar D_i(a) + \ldots , \quad a \in M_N(C^\infty(U))$$is an isomorhism of algebras.
\end{definition}

We  want to do the same construction for deformation quantizations
as we did for endomorphism bundles. We therefore need the infinitesimal
version of $M_N(\cW_n)$. This is just given by considering
$\bbO_{2n}[[\hbar ]]$, $\bbO_{2n}=\bbC [[\hat{x}_1,\ldots,\hat{x}_n ,
\hat{\xi}_1,\ldots, \hat{\xi}_n ]]$,  with the same product as in the
Weyl quantization. With this product we denote the algebra by
$\bbA^\hbar$. The infinitesimal structure of $ M_N(\cW_n)$ will then be
$M_N(\bbA^\hbar)$.

\begin{definition}
A formal symplectomorphism of $\bbO_{2n}$ is a continuous
automorphism of $\bbO_{2n}$ induced from an automorphism of
$$\bbR [[
\hat{x}_1 , \ldots , \hat{x}_n, \hat{\xi}_1 , \ldots , \hat{\xi_n}]]$$
that preserves the  formal standard Poisson bracket $\{ \cdot ,
\cdot \}$ on $\bbO_{2n}$.

Let $G$ be the subgroup of  automorphisms of $M_N(\bbA^\hbar )$ such that
$\Phi \in G$ if $\Phi$ is $\hbar$ linear and continuous. Moreover, $\Phi$ mod $\hbar$ becomes  an automorphism $\Phi_0$ of $M_N(\bbO_{2n})$ and $\Phi_0$
induces a formal symplectomorphism on $\bbO_{2n}$.
\end{definition}

If $\Phi \in G$ we let $\varphi $ denote the induced symplectomorphism
on $\bbO_{2n}$. In this case we will say that $\Phi$ is a
automorphism over $\varphi$.

\begin{lemma}\label{inner} Every automorphism of $M_N(\bbA^\hbar )$ over the identity symplectomorphism is inner.
\end{lemma}

\textit{Proof.} Let $\Phi$ be such an automorphism. Since it is an automorphism mod $\hbar$ over the identity, it is inner mod $\hbar$, and we can hence assume that $\Phi$ is the identity mod $\hbar$. In other words,
$$\Phi (a)= a+\hbar D_1(a) + \ldots, \quad a \in M_N(\bbO_{2n}). $$
Since $\Phi$ is an automorphism, $D_1$ is a derivation of
$M_N(\bbO_{2n})$ and hence of the form $X+[A, \cdot ]$, where $X$ is a
formal vector field and $A \in M_N(\bbO_{2n})$. If we assume that
$a,b\in \bbO_{2n}$,  we have that
$$\{ a, D_1(b)\} + \{ D_1(a),b\} = D_1 (\{a,b \})$$
 since $\Phi$ is an automorphism. This means that $X$ is a formal hamiltonian vector field. Therefore there exists an element $x$ in $\bbO_{2n}$ such that
 $D_1(a)=\{ x,a \}$. Hence we have
$$\Phi \circ \exp (-ad (x+A))=id \hbox{ mod } \hbar^2.$$

Continuing in this way, the result follows.
\\

It is well known from \cite{NT1}, that the Lie algebra  of $G$,
in the case where $N=1$, is given by
$$
\fg_0=\left\{ \frac{ia}{\hbar}| a \in \bbA^\hbar , \quad a \hbox{ real mod } \hbar, \quad a \in (\hat{x}_1, \ldots , \hat{\xi}_1 , \ldots )^2 \hbox{ mod } \hbar \right\}/\frac{\bbC[[\hbar ]]}{\hbar}
$$
and that any element of $G$ is of the form $\exp ( g)$,  $g \in\fg_0$.

We therefore see that the Lie algebra $\fg_0$ of $G$ for arbitrary $N$ is given by
\begin{eqnarray*}
\fg_0 &=& \left\{ \frac{ia}{\hbar} + b | b \in M_N(\bbA^\hbar ), \quad a \in \bbA^\hbar , \quad a \hbox{ real mod } \hbar ,  \right. \\
&& \left.  a \in (\hat{x}_1, \ldots , \hat{\xi}_1 , \ldots )^2 \hbox{ mod } \hbar \right\}/\frac{\bbC[[\hbar ]]}{\hbar}.
\end{eqnarray*}
To this Lie algebra we add the derivations $\partial_{\hat{x}_1}, \ldots, \partial_{\hat{\xi}_1}, \ldots $ and call the enlarged Lie algebra $\fg$.

Let us  suppose that we are given a deformation quantization $A^\hbar_E$ of $End(E)$. We define
$$I_m= \{ a \in A^\hbar_E |a(m)=0 \},$$
and we let $I_m^n$ denote the $n$'th power of the ideal $I_m$ in the undeformed product. The jet of $A^\hbar_E$ in $m$ is defined by
$$ J_mA^\hbar_E = \lim_{\stackrel{\leftarrow}{k}} A^\hbar_E / I_m^k.$$

Since the value of the product in $A^\hbar_E$ in a point only depends
on the derivatives in that point, the product descends to $J_mA^\hbar_E$.

If we choose a trivialization of $A^\hbar_E$ around $m$, we get an
isomorphism $J_mA^\hbar_E \isomoto M_N(\bbA^\hbar ) $. Any other
trivialization will give an automorphism of $M_N(\bbA^\hbar )$ in $G$.

As in the case of $End(E)$ we do the following
\begin{definition} The frame bundle $\tilde{M}_{A^\hbar_E}$ is given by
$$\tilde{M}_{A^\hbar_E} =\{ (m,\Phi )| m \in M, \Phi: M_N(\bbA^\hbar ) \isomoto J_m\bbA^\hbar_E \}.$$
\end{definition}

As before, $\tilde{M}_{A^\hbar_E}$ is a principal bundle with fibre  $G$.

\begin{proposition}For all $(m,\Phi ) \in \tilde{M}_{A^\hbar_E} $ there exists an isomorphism
$$\omega : T_{(m,\Phi )} \tilde{M}_{A^\hbar_E} \rightarrow \fg $$
 satisfying
\begin{eqnarray*}
(i) \quad &  \omega (A^* )=A & A \in \fg_0 \\
(ii) \quad & \varphi^* \omega = ad \varphi^{-1} \omega & \varphi \in G\\
(iii) \quad& d \omega + \frac{1}{2} [\omega , \omega ]=0
\end{eqnarray*}
where $A$ is the fundamental vector field corresponding to $A$.

In other words $\omega$ is a flat connection with values in $\fg$.
\end{proposition}

\textit{Proof.} The same as  the case of endomorphism bundles.
\\

As in the case of endomorphism bundles, we get a flat connection $\nabla$ in the bundle
$$JA^\hbar_E =\tilde{M}_{A^\hbar_E} \times_G M_N(\bbA^\hbar ).$$
\begin{proposition} The complex $(\Omega^* (M, JA^\hbar_E ), \nabla )$ is acyclic and $$\ker \nabla \isomoto A^\hbar_E.$$
\end{proposition}

\textit{Proof.} The same  as  the case of endomorphism bundles.
\\

We see that $H=U(n)\times SU(N)$ is a maximal compact subgroup of $G$. The $U(n)$-component comes from a maximal compact
subgroup of  the symplectic group, $SP(2n)$, and  $SU(N)$ comes from the maximal compact
subgroup of the action of $Gl_N$ on $M_N(\bbC )$.

Since $H$ is a maximal compact subgroup of
$G$, we can reduce the bundle $\tilde{M}_{A^\hbar}$
to an $H$ bundle $P$, which is easily seen to be a
reduction of the principal bundle consisting of dual symplectic frames
and frames of $End (E)$. We thus see that $JA^\hbar_E$ is, in fact,
isomorphic to $P\times_H\bbA^\hbar$. We will
denote this bundle by $\bbW$.

We introduce a grading on $\bbA^\hbar$ in which $\hat{x}_1, \ldots ,
\hat{x}_n, \hat{\xi}_1, \ldots , \hat{\xi}_n$ has degree 1 and $\hbar$
has degree 2. This also gives a grading on $M_N(\bbA^\hbar
)$. Furthermore we see that the action of $U(n)$ on  $\bbA^\hbar$
preserves the grading and we hence get a grading on $\bbW$.

We note that we have an extension of Lie algebras
\begin{equation}\label{ext} 0 \rightarrow \frac{1}{\hbar}i\bbR +\bbC
  [[\hbar ]] \rightarrow
\tilde{\fg} \rightarrow \fg \rightarrow 0 \end{equation}
where
$$
\tilde{\fg} = \left\{ \frac{ia}{\hbar} + b \mid a \in \bbA^\hbar, \quad
  a \hbox{ real mod }\hbar, \quad b \in M_N(\bbA^\hbar ) \right\},
$$
with bracket given by commutators. The grading $M_N(\bbA^\hbar)$ also gives a grading on $\tilde{\fg}$ and therefore also on the subbundle of $\bbW/\hbar$ with fibers $\tilde{\fg}$, denoted by $\tilde{\fg}_M$. Let $\fh = \fu (n) + \fsu (N)$. There is an embedding of $\fh $ in $\tilde{\fg}$ compatible with the quotient map $\tilde{\fg}\rightarrow \fg$ and the embedding  $\fh \rightarrow \fg$. This embedding is obvious in the case of $\fsu (N)$. For the case $\fu (n)$, this embeds in ${\frak{sp}}(2n)$, the Lie algebra of the symplectic group. Moreover ${\frak{sp}}(2n)$ embeds in $\fg$ as the Lie sub algebra generated by elements on the form
$$i\hbar^{-1}v_1v_2, \quad v_1,v_2 \in \{ \hat{x_1}, \ldots , \hat{\xi}_1,\ldots \}.$$   The embedding of $\fu (n)$ in $\tilde{\fg}$ is  the restriction of the embedding of ${\frak{sp}}(2n)$ in $\tilde{\fg}$ given by 
$$ i \hbar^{-1}v_1v_2 \rightarrow i \hbar^{-1}v_1 *v_2+\frac{i}{2}\omega_{st}(v_1,v_2), \quad   v_1,v_2 \in \{ \hat{x_1}, \ldots , \hat{\xi}_1,\ldots \}.$$

Since $H$ acts
semi-simple on $\fg$, we get an $H$-equivariant lift of
the quotient map $\tilde{\fg}
\rightarrow \fg$. We can therefore lift the connection $\nabla$
to a connection $\tilde{\nabla}$ taking values in $\tilde{\fg}$.
This means that we have a collection of local $\tilde{\fg}$-one forms $\{A_i \}$, $i$ being
labels of trivializations of $\bbW$, satisfying
$$A_i=g_{ij}dg_{ji}+g_{ij}A_j g_{ji}$$
where $g_{ij}$ are the transition functions.

This
connection however, is not flat. But because of the extension
\ref{ext}, the
curvature $dA_i+\frac{1}{2}[A_i,A_i]$ is in $ \Omega^2(M,\frac{1}{\hbar}i\bbR +\bbC [[\hbar
]])$. Clearly the associated cohomology class is independent of the
choice of lifting of $\nabla$. By checking the definition
of $\nabla$, one sees that the $\frac{1}{\hbar}$ component of
$\tilde{\nabla}^2$ is $\frac{\omega}{i\hbar}$, where $\omega$ is the
symplectic structure.

 So for each deformation quantization of $End(E)$ we have a connection
 $\nabla$ in $\bbW$, such that $\ker \nabla$ is isomorphic to the
 deformation quantization. And the lifting of $\nabla$ to a
 $\tilde{\fg}$-valued connection gives an element in
$$ \frac{\omega}{i\hbar}+H^2(M,\bbC [[\hbar ]]).$$

\begin{proposition} \label{con}
Let $A^\hbar_{1,E}$ and $A^\hbar_{2,E}$ be
  deformation quantizations with characteristic classes $\theta_1$ and
  $\theta_2$. Then $A^\hbar_{1,E} \simeq A^\hbar_{2,E}$ if and only if
  $\theta_1=\theta_2$ in $\frac{\omega }{i\hbar}+H^2(M,\bbC [[\hbar
  ]])$.
\end{proposition}

The proof of the above proposition relies on the following

\begin{lemma} Let $I_\omega : TM \rightarrow T^*M$ be the bundle
isomorphism induced by $\omega$. Since $T^*M \subset \bbW$, the isomorphism $I_\omega
$ induces an element $A$ in $\Omega^1(M,\bbW)$. Put
$$A_{-1} =\frac{A}{\hbar} \in \Omega^1 (M,\tilde{\fg}_M ).$$
Then the complex
$$(\Omega^* (M,\tilde{\fg}_M ), Ad A_{-1} )$$ is acyclic.
\label{kozul}
\end{lemma}

\textit{Proof.} Since the action of $ Ad A_{-1} $ commutes with the action of $C^\infty (M)$ on $\Omega^* (M,\tilde{\fg}_M )$ it is enough to prove the statement locally. Locally the complex is just 
$$\left( C^\infty (\bbR^{2n}) \otimes (i/\hbar \hat{\Omega}^* +M_N(\hat{\Omega}^*)[[\hbar ]]),\id \otimes \hat{d}\right),$$
where $(\hat{\Omega^*}, \hat{d})$ is the the complex of formal differential forms. From this the lemma follows.
\\
%From Lemma 3.12 in \cite{NT4} we know, in  the case of
%a trivial bundle $E$, that
% $(\Omega^*(M, \bbW),Ad(A_{-1}))$ is acyclic and the cohomology is
%isomorphic to $C^\infty (M)[[\hbar ]]$. Since
%\begin{eqnarray*}
%&&(\Omega^*(M,\bbW), Ad(A_{-1}))= \\
%&&(\Omega^*(M, \bbW),
%Ad(A_{-1}))\otimes_{C^\infty (M)[[\hbar ]]} \Gamma (End(E))[[\hbar ]]
%\end{eqnarray*}
%and $\Gamma (End(E))[[\hbar ]]$ is projective over $C^\infty (M)[[\hbar ]]$, the
%result follows.

\textit{Proof(of Proposition).} Let $\tilde{\nabla}_1$ and $\tilde{\nabla}_2$ be the
two connections with
$\tilde{\nabla}_1^2=\tilde{\nabla}_2^2=\theta$. Note that we can
actually assume that $\tilde{\nabla}_1^2=\tilde{\nabla}_2^2$ in
$\frac{ \omega}{i\hbar}+\Omega^2(M,\bbC [[\hbar ]])$. We have
$$ \tilde{\nabla}_1-\tilde{\nabla}_2=R_0+R_1+ \ldots , \quad R_i \in
\Omega^1(M,\tilde{\fg}_M^i )$$
From the equality of the curvatures we get $[A_{-1}, R_0]=0$ and hence
by  lemma \ref{kozul} we have an element $g_1 \in  \Gamma (\tilde{\fg}_M^1)$, such that $[g_1,A_{-1}]=R_0$.
Therefore considering the connection $\nabla_{2,0}=ad(\exp
g_1)\nabla_2$ we have
$$\tilde{\nabla}_1-\tilde{\nabla}_{2,0}=R_1'+\ldots, \quad R_i'\in
\Omega^1(M,\tilde{\fg}_M^i )$$
Continuing by induction and using the Hausdorff-Campbell formula we get an element $g$ in $\Gamma (Aut(\tilde{\fg}_M))$
conjugating $\nabla_1$ into $\nabla_2$, and hence the deformations
will be isomorphic.

Let us now assume that  $A^\hbar_{1,E}$ and $A^\hbar_{2,E}$ are
isomorphic. This induces an isomorphism between
$\tilde{M}_{A^\hbar_{1,E}}$ and $\tilde{M}_{A^\hbar_{2,E}}$ compatible with the
connections on $\tilde{M}_{A^\hbar_{1,E}}$ and
$\tilde{M}_{A^\hbar_{2,E}}$ . In particular we get an automorphism of
$\bbW$ mapping $\nabla_1$ to $\nabla_2$, from which we see that
$A^\hbar_{1,E}$ and $A^\hbar_{2,E}$ have the same characteristic class.

\begin{thm} The deformation quantizations of an endomorphism
  bundle are classified by the affine space
$$\frac{\omega}{i\hbar} + H^2(M,\bbC [[\hbar ]]).$$
\end{thm}

\textit{Proof.} We only need to prove that for a class $\theta$ in
$\frac{ \omega}{i\hbar} +\Omega^2 (M, \bbC [[\hbar ]])$, we have a
deformation quantization with characteristic class $\theta$. To do
this we start with a connection in $P$ and thus get an $H$-connection $\nabla$ in $\bbW$. We have that
$$[\nabla + A_{-1},\nabla + A_{-1}]=\frac{ \omega}{i \hbar}+2 [\nabla ,
A_{-1}]+[\nabla , \nabla ].$$

One checks that $[A_{-1}, [\nabla , A_{-1}]]=0$, and according to lemma
\ref{kozul} we get an element $A_0 \in \Omega^1 (M,\tilde{\fg}_M )$ such that $[A_{-1}, A_0] =
[\nabla , A_{-1} ]$. If we put $\nabla_0=\nabla+A_{-1}+A_0$ we have
$$[\nabla_0 , \nabla_0 ]- \theta= 0 \hbox{ mod }
\Omega^2(M, \tilde{\fg}_M^{i\geq 0}).$$

Let us now assume, that we have constructed $\nabla_n$ with
$\nabla_n^2-\theta=0$ mod $\Omega^2 (M, \tilde{\fg}_M^{i\geq n})$. By the Bianchi
identity we have $[ \nabla_n
, [\nabla_n, \nabla_n ]]=0$, and therefore $[A_{-1},([\nabla_n,
\nabla_n]-\theta)_n ]=0$, where $([\nabla_n , \nabla_n ]-\theta )_n$
is the $n$-th component of $[\nabla_n ,  \nabla_n ] -\theta$. As
before we get an element $A_n$ with $[A_{-1},A_n]=([\nabla_n,
\nabla_n]-\theta)_n$, and considering $\nabla_{n+1}=\nabla_n+A_n$ we
have $[\nabla_{n+1},\nabla_{n+1}]-\theta=0$ mod $\Omega^2(M,
\tilde{\fg}_M^{i\geq n+1})$.

We can thus for each class $\theta$ in $\frac{\omega}{i \hbar }
+\Omega^2(M,\bbC[[\hbar]])$ construct a connection $\nabla_F$ with values in
$\tilde{\fg}$ and curvature $\theta$. We therefore only need to check
that the complex $(\Omega^*(M,\bbW), \nabla_F)$ is acyclic, that the
kernel is isomorphic to $\Gamma (End(E)[[\hbar]]$ and that the product
induced by the product on $\bbW$ gives a deformation quantization of $End(E)$.
These results, however follow from the proof of the Proposition \ref{con},
since we locally can conjugate $\nabla_F$ to a connection on the form
$$d-\sum_{i=1}^n ( \partial_{\tilde{x}_i}\otimes dx_i +
\partial_{\tilde{\xi}_i}\otimes d \xi_i ).$$

\section{Lie Algebra Cohomology} \label{liesekt}
\begin{definition}
A differential graded Lie algebra $(\fg , d )$ over a commutative
unital ring $k$ is a ($\bbZ / 2 \bbZ $
or $\bbZ$) graded $k$-module $\fg$ with a bracket operation $[\cdot,
\cdot ] : \fg^j \times \fg^i \rightarrow \fg^{i+j}$ and a differential
$\partial : \fg^i \rightarrow \fg^{i-1}$ satisfying :
\begin{eqnarray*}
(i) &\partial [g_1,g_2 ]=[\partial g_1,g_2 ] + (-1)^{|g_1|}[g_1,\partial g_2], \\
(ii)& [g_1,g_2]=-(-1)^{|g_1||g_2| }[g_2,g_1], \\
(iii)&[g_1,[g_2,g_3]]+(-1)^{|g_3|(|g_1|+|g_2|)}[g_3, [g_1,g_2]]\\
&+(-1)^{|g_1|(|g_2|+|g_3|)}[g_2,[g_3,g_1]]=0,
\end{eqnarray*}
where $|\cdot |$ is the degree.
\end{definition}

A $\fg$ module $\bbL^*$ is a complex $\bbL^*$ with an action of $\fg$,
i.e. we have a map $\fg^i \times \bbL^j \rightarrow
\bbL^{i+j}$ satisfying
$$g_1g_2l-(-1)^{|g_1||g_2|}g_2g_1l=[g_1,g_2]l$$
and
$$\partial_{\bbL^*}
(gl)=(\partial_{\fg}g)l+(-1)^{|g|}g(\partial_{\bbL^*}l).$$

Given a differential graded Lie algebra $\fg$, we can define a
differential graded Lie algebra $\fg [\epsilon ]$ as follows
\begin{eqnarray*}
\bullet & \fg [\epsilon ] = \fg + \epsilon \fg & \hbox{ where }|\epsilon | =1, \\
\bullet & [g_1, \epsilon g_2 ] = \epsilon [g_1, g_2], \\
\bullet & [\epsilon g_1, \epsilon g_2 ]=0,\\
\bullet & \partial (g_1+\epsilon g_2 )=\partial_{\fg}g_1+g_2 -\epsilon \partial_\fg g_2.
\end{eqnarray*}

Also one can construct the enveloping algebra  by setting
$$U(\fg ) = T(\fg)/ (g_1\otimes g_2 -(-1)^{|g_1||g_2|}g_2 \otimes g_1
=[g_1, g_2])$$
where $T(\fg)$ is the tensor algebra. Furthermore $U(\fg)$ has a
differential induced by the differential on $\fg$ by the graded
Leibniz rule and a grading.

We note that $(U(\fg[\epsilon]), \partial)$ is a $\fg$ module. Let $\fh$
denote a Lie sub algebra of $\fg$. For a $\fg$
module $\bbL^*$  we define
$$C^*(\fg, \fh ; \bbL^*)=\Hom_{U(\fg)}(U(\fg [\epsilon ]) \otimes_{U(\fh
  [\epsilon ])}k, \bbL^*).$$
This will have a differential induced by the differential on
$U(\fg[\epsilon])$ and the differential on $\bbL^*$. The homology
of this complex will be denoted by $H^*(\fg, \fh ; \bbL^*)$.

We are now going to give a construction of  classes in $C^*(\fg, \fh
; \bbL^*)$ in special cases.
%The Weil algebra of $\fh$ is defined by
%$$W^*(\fh )= C^*(\fh [\epsilon ], \fh ; k)$$
%We can regard $W^*(\fg )$ as being the free graded commutative algebra
%generated by two copies of $\fh'$, one copy $\fh'$ in degree one and
%one copy $\epsilon\fh'$ in degree
%two.
First we assume that $\bbL^*$ is  homotopically constant in
the sense defined below
\begin{definition}
A $\fg$-module $\bbL^*$ is called  homotopically constant if there exist
operations
\begin{eqnarray*}
& \iota_g : \bbL^* \rightarrow \bbL^{*-1} & g\in \fg 
\end{eqnarray*}
satisfying
\begin{eqnarray*}
& [\partial , \iota_g ] = L_g; \quad [\partial , L_g]=0: \\&  [L_{g_1}, \iota_{g_2}
]=\iota_{[g_1,g_2]}; \quad [\iota_{g_1}, \iota_{g_2} ]=0:
\end{eqnarray*}
where we have denoted the action of $\fg$ by $L_g$.

In other words, we have an action of the differential graded algebra
$U(\fg[ \epsilon ])$ on $\bbL^*$.
\end{definition}

%We thus get a morphism of complexes $\bbL^* \rightarrow C^*(\fg ,
%\bbL^*)$ given by
%$$l\rightarrow \sum \lambda_nl$$
%where
%$$\lambda_nl (g_1 \wedge \ldots  \wedge g_n) = \iota_{g_{n}}\cdots
%\iota_{g_{n}}l$$

If we furthermore assume that there is an $\fh$-equivariant projection
$\nabla: \fg \rightarrow \fh$ of the embedding $\fh \rightarrow \fg$,
we get the usual Chern-Weil homomorphism, i.e. a map of complexes
$$CW: C^*(\fh[\epsilon],\fh;\bbL^* ) \rightarrow C^*(\fg,\fh;\bbL^* )$$
given in the following way:

For an element $g_1\wedge g_2$ define  $R(g_1\wedge g_2)=[\nabla
(g_1), \nabla (g_2)]-\nabla ([g_1,g_2])$. Taking cup product gives
$R^n:\wedge^{2n}\fg \rightarrow \wedge^n\epsilon \fh$. By composition this
gives a map
$$\varphi:C^*(\fh[\epsilon],\fh;\bbL^*) \rightarrow C^*(\fg,\bbL^*).$$
There are operations $\lambda^n$ on $\bbL^*$ given by
$$g_1\wedge \ldots \wedge g_n \rightarrow \iota_{g_1-\nabla g_1}\cdots
\iota_{g_n-\nabla g_n} l,$$
where $l \in \bbL^*$. Finally we set
$$CW( a)=\sum_n \lambda^n\cup \varphi (a )$$
where $\cup$ is the cup product. One checks that this gives a morphism
of complexes.

%\begin{eqnarray*}
% \fh' \ni \varphi & \rightarrow &( g \rightarrow  \varphi (\nabla g)) \\
% \epsilon \fh' \ni \varphi &\rightarrow & (g_1\wedge g_2 \rightarrow
% \varphi (R(g_1,g_2)))
%\end{eqnarray*}
% where $R(g_1,g_2)=[\nabla (g_1), \nabla (g_2)]-\nabla ([g_1,g_2])$.

%Using the above morphism $\bbL^* \rightarrow C^*(\fg ; \bbL^*)$, the
%Chern-Weil homomorphism and the cup product we get a morphism of
%complexes
%$$CW : W^*(\fh) \otimes \bbL^* \rightarrow C^*(\fg; \bbL^*)$$

%We note that we have operations $\iota_h,L_h$ on the two complexes
%$W^*(\fg) \otimes \bbL^*$ and $C^*(\fg; \bbL^*)$ and that $CW$
%preserves these operations. Elements in $W^*(\fg) \otimes \bbL^*$ for
%which the operations are zero is denoted $(W^*(\fg) \otimes
%\bbL^*)^{basic}$ and elements in  $C^*(\fg; \bbL^*)$ that are zero
%under the operations are easily seen to be  $C^*(\fg, \fh; \bbL^*)$

Next we will give a construction of classes in $C^*(\fh [\epsilon ],
\fh ;\bbL^*)$ for special cases of $\bbL^*$. To this end we need the following
\begin{definition}
A $\fg$-module $\bbL^*$ is called  very homotopically constant
if $\bbL^*$ is  homotopically constant and we have
operations
\begin{eqnarray*}
& L_{\underline{g}}:\bbL^* \rightarrow \bbL^{*+1} & g \in \fg \\
& \iota_{\underline{g}}:\bbL^* \rightarrow \bbL^* & g \in \fg
\end{eqnarray*}
satisfying the eight conditions
\begin{eqnarray*}
& [\partial , \iota_{\underline{g}} ]= L_{\underline{g}} - \iota_g;  \quad
[\partial, L_{\underline{g}}]=L_g; \quad [\iota_{\underline{g}_1}, \iota_{\underline{g}_2}]=0;
\quad   [\iota_{g_1}, \iota_{\underline{g}_2}]=0; \\ &
[\iota_{\underline{g}_1}, L_{g_2}]=\iota_{\underline{[g_1,g_2]}}; \quad
[\iota_{\underline{g}_1}, L_{\underline{g}_2}]=0; \quad [L_{g_1},
L_{\underline{g}_2} ]=L_{\underline{[g_1,g_2]}}; \quad
  [L_{\underline{g}_2} , L_{\underline{g}_2} ]=0;
\end{eqnarray*}
In other words, we have an action of the differential graded algebra
$U(\fg[\epsilon, \eta ])=U(\fg[\epsilon][\eta ])$ on $\bbL^*$.
\end{definition}
We now assume that $\bbL^*$ is a very homotopically constant
$\fh$-module. We denote by
$\bbL^*_{\fh +\underline{\fh}}$ the elements $l \in \bbL^*$ with $L_hl=0$ and
$L_{\underline{h}}l=0$. We note that this is a complex. We get a morphism
of complexes $\bbL^*_{\fh +\underline{\fh}} \rightarrow
C^*(\fh[\epsilon ] , \fh ; \bbL^*)$
by:
$$\bbL^*_{\fh +\underline{\fh}} \ni l \rightarrow (h_1\cdots  h_n
\rightarrow \iota_{\underline{h}_1} \cdots \iota_{\underline{h}_n} l).$$

\subsection{Examples}

We are  going to give some examples of relative classes in Lie
algebra cohomology.

\textbf{Example 1} Consider the extension
\ref{ext} and choose an $\fh$-equivariant lift $\nabla : \fg
\rightarrow \tilde{\fg}$ of the quotient map $\tilde{\fg}\rightarrow
\fg$, where $\fh$ is ${\frak{u}}(n) + {\frak{su}}(N)$. We then define the
class $\theta$ in $C^*(\fg,\fh;\frac{1}{\hbar}\bbC [[ \hbar ]] )$ by
\begin{equation}
\theta (g_1,g_2)=[\nabla g_1, \nabla g_2]-\nabla([g_1,g_2]). 
\end{equation}
We want to show that $\theta$ actually comes from a sort of
Chern-Weil map. Let $k:\tilde{\fg} \rightarrow \fg$ denote the
quotient map and let $\tilde{\fh}=k^{-1}(\fh)$, i.e.
$$\tilde{\fh}=\fh+\frac{i}{\hbar}\bbR+\bbC [[ \hbar ]]$$
Note that
$C^*(\fg,\fh; \bbL^*)$ and $C^*( \tilde{\fg}, \tilde{\fh}; \bbL^*)$
  are quasi-isomorphic when $\bbL^*$ is a $\fg$ module. We have a
  Chern-Weil homomorphism
$$CW: C^*( \tilde{\fh}[\epsilon ], \tilde{\fh}; \bbL^*) \rightarrow
C^*(\tilde{\fg}, \tilde{\fh}; \bbL^*)$$
as before. A choice of an $\tilde{\fh}$ equivariant split $\nabla ' :
\tilde{\fg} \rightarrow \tilde{\fh}$ is given by $\nabla '=\nabla ''
\circ k + id -\nabla \circ k$, where $\nabla'': \fg \rightarrow \fh$
is an $\fh$-equivariant splitting of the embedding $\fh \rightarrow
\fg$. Let $\theta$ be   
 the projection of  $\tilde{\fh}$ on $\frac{i}{\hbar}\bbR+\bbC [[ \hbar ]]$. It is now
easy to see that  under the quasi isomorphism
between $C^*(\fg,\fh; \frac{1}{\hbar}\bbC[[\hbar]])$ and $C^*( \tilde{\fg}, \tilde{\fh};
\frac{1}{\hbar}\bbC [[\hbar ]])$, the class $CW(-\theta )$ is the same as the $\theta$ we defined in the start of this example.

\textbf{Example 2} Some other classes  in  $ C^*(\fg ,\fh)$ we need are also coming
from a Chern-Weil construction.  We  consider an $\fh$-equivariant splitting $\nabla' $ of the embedding
 $\fh \rightarrow \fg$.
Composed with the the projection
 ${\frak{u}}(n)+{\frak{su}}(N) \rightarrow {\frak{su}}(N)$ we get an
 $\fh$ equivariant map $\nabla : \fg \rightarrow {\frak{su}}(N)$. Using
 this we therefore get a Chern-Weil homomorphism
$$CW: C^*(\fsu (N)[\epsilon],\fsu(N))  \rightarrow C^*(\fg,\fsu (N)).$$
It is clear that this homomorphism in fact maps into $C^*(\fg,
\fh )$. We therefore get the usual classes, for example the usual
chern character $ch$, which is
$$\exp (R)=\sum \frac{1}{n!}Tr(R^n),$$
where
$$R(g_1,g_2)=[\nabla g_1, \nabla g_2 ]-\nabla[g_1,g_2].$$
Here $Tr$ denotes the usual  normalized trace on $\fsu (N)$.
This class is of course the Chern Weil map on the symmetric
polynomium $ch$ on $\fsu (N)$ given by
$$ch(h_1, \ldots , h_n)=\frac{1}{n!}\sum_{\sigma \in
    S_n}\frac{1}{n!}Tr(h_{\sigma (1)} \cdots h_{\sigma (n)}).$$

Since $C^*(\fsu (N)[\epsilon ] , \fsu (N))$ embeds in
$C^*(\tilde{\fh}[\epsilon ], \tilde{\fh})$, the Chern Weil construction
given in this example is just a particular case of the Chern Weil
homomorphism
$$CW: C^*( \tilde{\fh}[\epsilon ], \tilde{\fh }) \rightarrow
C^*(\tilde{\fg}, \tilde{\fh})$$

\textbf{Example 3} We can of course do the construction from
example 2 for $\fu (n)$ instead of
$\fsu (N)$, and therefore get a Chern-Weil homomorphism
$$C^*(\fu (n)[\epsilon],\fu (n))\rightarrow C^*(\fg ,\fh).$$
This again can also be viewed as the composition
$$C^*(\fu (n)[\epsilon],\fu (n))\rightarrow C^*(\tilde{\fh}[\epsilon ]
,\tilde{\fh})\rightarrow C^*(\tilde{\fg} ,\tilde{\fh})$$
We will in particular be interested in the symmetric polynomium $\hat{A}$ coming from the map
$$h \rightarrow \det \left( \frac{h/2}{\sinh (h/2)}\right)$$

\section{Cyclic Homology}\label{cyclic}
We consider a differential graded unital algebra $(A,\delta )$ over a
commutative ring $k$ containing $\mathbb{Q}$, i.e. an algebra $A$ that can be written as $A=\oplus_n A^n$, where $A^n$'s are independent $k$-submodules of $A$ and $A^nA^m \subset A^{n+m}$. Elements in $A^n$ is said to have degree $n$ and we will denote the degree of an element $a$ by $|a|$. Furthermore $\delta
: A^* \rightarrow A^{*-1}$ has to be a differential  and satisfy $\delta (ab)=\delta (a)b+
(-1)^{|a|}a\delta (b)$.

Define an operator $\tau$ on $A^{\otimes (n+1)}$ by
$$\tau (a_0 \otimes \ldots \otimes
a_n)=(-1)^{(|a_n|+1)\sum_{i=0}^{n-1}(|a_i|+1)}a_n \otimes a_0 \otimes
\ldots \otimes a_{n-1},$$
and consider the complex
$$\ldots \stackrel{b + \delta}{\longleftarrow}
C^n(A)/Im(1-\tau)\stackrel{b+\delta}{\longleftarrow}
C^{n+1}(A)/Im(1-\tau)\stackrel{b+\delta}{\longleftarrow} \ldots $$
where $C^{n+1}(A)$ denotes the set of  elements of the form $a_0 \otimes \ldots \otimes
a_k$ with $k+\sum |a_i|=n$,
\begin{eqnarray*}
b(a_0 \otimes \ldots \otimes a_n)=\sum_{k=0}^{n-1}(-1)^{k+\sum_{i=0}^k
  |a_i|}a_0 \otimes \ldots \otimes a_ka_{k+1}\otimes \ldots \otimes
a_n \\
+ (-1)^{(|a_n|+1)\sum_{i<n}(|a_i|+1)+|a_n|}a_na_0\otimes \ldots
\otimes a_{n-1}
\end{eqnarray*}
and
$$\delta (a_0 \otimes \ldots \otimes
a_n)=\sum_{k=0}^{n}(-1)^{\sum_{i=1}^{k-1}(|a_i|+1)} a_0\otimes \ldots \otimes
  \delta(a_k) \otimes \ldots \otimes a_n$$
The complex is denoted by $C^\lambda_*(A)$, and the homology is the cyclic
homology of $A$ denoted by $HC_*(A)$.

The reduced cyclic homology is given by the homology of the complex
$$\ldots \stackrel{b + \delta}{\longleftarrow}
\overline{C}^n(A)/Im(1-\tau)\stackrel{b+\delta}{\longleftarrow}
\overline{C}^{n+1}(A)/Im(1-\tau)\stackrel{b+\delta}{\longleftarrow} \ldots $$
where $\overline{C}^*(A)$ comes from considering $\bar{A}^{\otimes *}$
instead of $A^{\otimes *}$, where $\bar{A}=A/k \cdot 1$.

The reduced cyclic homology of $A$ is denoted by $\overline{HC}_*(A)$, and the  complex above, computing the reduced cyclic homology, is denoted by $\overline{C}^\lambda_*(A)$.

It is well known, see \cite{Lo} and \cite{BNT}, that there is an exact sequence
\begin{equation} \label{lang}
\ldots HC_n(k) \rightarrow HC_n(A) \rightarrow \overline{HC}_n(A)
\rightarrow
HC_{n-1}(k) \rightarrow \ldots
\end{equation}

We will briefly give a construction, due to Brodzki, of the connecting morphism  $\overline{HC}_*(A)
\rightarrow HC_{*-1}(k)$ at the level of complexes, see  \cite{Br}
and \cite{BNT}; i.e. a morphism
of complexes
$$Br: \overline{C}^\lambda_*(A) \rightarrow C^\lambda_{*-1}(k)$$
giving the connecting homomorphism at the level of homology. Let $l:A
\rightarrow k$ be a $k$-linear map with $l(1)=1$. Put
\begin{eqnarray*}
\rho (a) & = & l(\delta (a)), \qquad a \in A \\
\rho (a_1 \otimes a_2 ) &=& l(a_1)l(a_2)-l(a_1a_2), \qquad a_1 \otimes
a_2 \in A^{ \otimes 2}  \\
\rho & = & 0\quad  \hbox{ on }A^{\otimes m}\hbox{, } \quad m \geq 3
\end{eqnarray*}
and define $Br: \overline{C}^\lambda_*(A) \rightarrow
C^\lambda_{*-1}(k)$ by setting 
\begin{eqnarray*}
& Br(a_0 \otimes \ldots \otimes a_m)=
 \sum_{i=0}^m (-1)^{\sum_{k<i}(|a_k|+1)\sum_{k \geq
    i}(|a_k|+1)}
(\rho \otimes \ldots \otimes \rho) \\ & ( a_i \otimes \ldots \otimes a_0
\otimes a_m \otimes \ldots \otimes a_{i-1}) (n+1)!\cdot 1^{\otimes 2n+1}
\end{eqnarray*}
on $\overline{C}^\lambda_{2n+1}$ and letting $Br$ be zero on 
  $\overline{C}^\lambda_{2n}$ $Br$.

We now consider the differential graded algebra $k[\eta]$, where
$\eta$ has degree one, $\eta^2=0$ and the differential is given by $\partial_\eta$.
For a differential graded algebra $A$ we define $A[\eta]$ to be $A
\otimes k[\eta]$, where $\otimes$ is the tensor product of
differential graded algebras. It is not difficult to see that
$HC_*(A[\eta ])=0$, and we therefore have
\begin{proposition}
The morphism

$$Br: \overline{C}^\lambda_*(A[\eta]) \rightarrow C^\lambda_{*-1}(k)$$
is a quasi isomorphism.
\end{proposition}

Since it is not standard we mention, that the reduced cyclic homology
is Morita invariant, at least in the case of algebras and matrices
over these algebras. To see this, let $A$ be an algebra, and let $l: A
\rightarrow k$ be a map needed in the
construction of $Br$. Let $tr$ denote the normalized trace $tr:M_n(A)
\rightarrow A$. We now have a commutative diagram
$$
\begin{array}{ccccccc}
 \rightarrow & HC^*(M_N(A)) & \rightarrow &
\overline{HC}^*(M_N(A)) & \stackrel{Br}{\rightarrow} & HC^{*-1}(k)&
\rightarrow \\
& \downarrow tr & & \downarrow tr && \big{\|} &
\\ \rightarrow & HC^*(A) & \rightarrow &
\overline{HC}^*(A) & \stackrel{Br}{\rightarrow} & HC^{*-1}(k)&
\rightarrow
\end{array}$$
where $Br:\overline{HC}^*(M_N(A)) \rightarrow HC^{*-1}(k)$ is induced
by $l \circ tr$. According to  \cite{Lo}, $tr:HC^*(M_n(A)) \rightarrow
HC^*(A)$ is an isomorphism. The result therefore follows from the exact sequence (\ref{lang}).

\subsection{Operations on the periodic complex}
For the periodic cyclic complex we consider $\prod_n A \otimes
\bar{A}^{\otimes n}$. We give this a $\bbZ/2 \bbZ$ grading by
$$|a_0 \otimes \ldots \otimes a_n |= n+\sum_i |a_i| \hbox{ mod }2.$$
On this we consider the differential $b+B+\delta$, where $b$ and
$\delta$ are given as before and
\begin{eqnarray*}
& B(a_0 \otimes \ldots \otimes a_n)\\
&=\sum_{i=0}^n (-1)^{\sum_{j\leq i}(|a_j|+1)\sum_{j\geq i+1}(|a_j|+1)} 1 \otimes a_i
\otimes \ldots \otimes a_n \otimes a_0 \otimes \ldots \otimes a_{i-1}.
\end{eqnarray*}
We will denote this complex by $CC^{per}_*(A)$.

The main feature about cyclic periodic homology that we are going to
need is the following;  see \cite{NT6} for complete formulas.
\begin{thm} \label{opera}
There is a morphism of complexes
$$\overline{C}^\lambda_*(A[\eta]) \otimes  CC^{per}_*(A) \rightarrow
CC^{per}_*(A)$$
satisfying the following:
\begin{itemize}
\item $n!\eta^{\otimes n+1} \cdot a= a$ for $ a \in CC^{per}_*(A)$
\item The component in $A$ of $(b_1 \otimes \ldots \otimes b_n) \cdot
  (a_0 \otimes \ldots \otimes a_m)$, where $b_1 \otimes \ldots \otimes
  b_n \in \overline{C}^\lambda_*(A)$, is zero when $m \not= n$ and equal to
$$\sum_i\frac{1}{n!}(-1)^{i(n-1)}a_0 [b_{i+1},a_1] \cdots [b_n, a_{n-i}][b_1,a_{n-i+1}]\cdots[b_i,a_n]$$
when $m=n$.
\end{itemize}

\end{thm}

The framework underlying Theorem \ref{opera} also gives other operations on
$CC^{per}_*(A)$, see \cite{NT6} for details. Let $(C^*(A,A), b)$ be the Hochshild cohomological
complex, i.e. $C^*(A,A)=Hom_k (\bar{A}^{\otimes *},A)$ and
\begin{eqnarray*}
& b \varphi (a_1, \otimes , a_{n+1})=(-1)^n a_1 \varphi(a_2, \ldots ,
a_{n+1})\\
&+ \sum_{j=1}^n (-1)^{n+j}\varphi (a_1, \ldots,a_j,a_{j+1}, \ldots ,
a_{n+1})-
\varphi(a_1, \ldots , a_n)a_{n+1}.
\end{eqnarray*}
Given two elements $\varphi$  in $C^n(A,A)$ and $\psi$ in $C^m(A,A)$, define
\begin{eqnarray*}
&\varphi \circ \psi (a_1, \ldots, a_{n+m-1}) \\
& =
\sum_{j \geq 0}(-1)^{(n-1)j}\varphi (a_1, \ldots , a_j , \psi
(a_{j+1}, \ldots , a_{j+m}), \ldots).
\end{eqnarray*}
Set
$$[\varphi, \psi ] = \varphi \circ \psi - (-1)^{(n+1)(m+1)} \psi \circ
\varphi.$$
With this bracket and with a suitably defined grading, $C^*(A,A)$
actually becomes a differential graded Lie algebra.

For $\varphi \in C^n (A,A)$ one can construct operations
\begin{eqnarray*}
& L_\varphi : CC^{per}_*(A) \rightarrow CC^{per}_{*-n+1} (A)\\
& I_\varphi : CC^{per}_*(A) \rightarrow CC^{per}_{*-n} (A)
\end{eqnarray*}
such that
\begin{eqnarray*}
& [ L_\varphi , L_\psi ]=L_{[\varphi , \psi ]},\\ 
&[ I_\varphi ,
L_\psi ]=I_{[\varphi, \psi ]}, \\ & [ B+b,I_\varphi]=I_{b \varphi}+L_\varphi .
\end{eqnarray*}

\section{The Fundamental Class}
We consider the reduced cyclic homology complex of
$M_N(\bbA^\hbar)[\hbar^{-1}]$. According to
Lemma 5.1.1 in \cite{BNT} and the Morita invariance of reduced the cyclic
homology, the homology is given in the following way
\begin{eqnarray*}
&\overline{HC}_i(M_N(\bbA^\hbar)[\hbar^{-1}])=\bbC[[\hbar, \hbar^{-1}],
\quad i=1,3, \ldots , 2n-1\\
&\overline{HC}_i(M_N(\bbA^\hbar)[\hbar^{-1}])=0 \quad \hbox{otherwise}
\end{eqnarray*}

A concrete generator for the homology in dimension $2n-1$ is given by
$$U_0=\frac{1}{2n(i\hbar)^n}\sum_{\sigma \in S_{2n}}(v_{\sigma
    (1)}\otimes \ldots \otimes v_{\sigma (2n)} ),$$
where $(v_1, \ldots , v_{2n})=(\hat{x}_1, \hat{\xi}_1, \ldots ,
\hat{x}_n, \hat{\xi}_n)$.

Let $\tilde{\fh}$ denote the inverse image of $\fh$ under the map
$\tilde{\fg} \rightarrow \fg$. Note that $U_0$ is invariant under the
action of $\tilde{\fh}$ and therefore, by the result on $\overline{HC}_*(M_N(\bbA^\hbar[\hbar^{-1}]))$, extends uniquely in homology to a class
$$ U\in C^*(\tilde{\fh}[\epsilon], \tilde{\fh};
\overline{C}^\lambda_*(M_N(\bbA^\hbar[\hbar^{-1}]))).$$
We wish to work with the differential graded algebra
$M_N(\bbA^\hbar[\hbar^{-1}])[\eta]$ instead of
$M_N(\bbA^\hbar[\hbar^{-1}])$. In the complex
$$\overline{C}^\lambda_*(M_N(\bbA^\hbar[\hbar^{-1}])[\eta])$$
we define $\eta^{(k)}= k!\eta^{\otimes k+1}$.

We also define operations on
$\overline{C}^\lambda_*(M_N(\bbA^\hbar[\hbar^{-1}])[\eta ])$  by 
\begin{eqnarray*} && \iota_g(a_0\otimes \ldots \otimes a_p)\\
&&= \sum_{i=0}^p(-1)^{\sum_{k \leq
    i}(|a_k|+1)(|g|+1)}
a_0 \otimes \ldots \otimes a_i \otimes g \otimes \ldots \otimes a_p,
\end{eqnarray*}
where $g \in M_N(\bbA^\hbar[\hbar^{-1}])[\eta ]$. Put $\iota_{\underline{g}}=\iota_{\eta g}$  and $L_{\underline{g}}=L_{\eta g}$ when  $g \in M_N(\bbA^\hbar[\hbar^{-1}])$. Here $L$ denotes the usual $L$ operation on the reduced cyclic complex. 

With these operations $\overline{C}^\lambda_*(M_N(\bbA^\hbar[\hbar^{-1}])[\eta ])$ becomes almost very homotopically constant over the commutator Lie algebra of  $ M_N(\bbA^\hbar[\hbar^{-1}])$ (Instead of the relations $[\partial , \iota_{\underline{g}} ]=L_{\underline{g}}-\iota_g$ and $[\partial , L_{\underline{g}} ]=L_g$ we have the relations  $[\partial , \iota_{\underline{g}} ]=L_{\underline{g}}+\iota_g$ and  $[\partial , L_{\underline{g}} ]=-L_g$. Phrased differently:  $\overline{C}^\lambda_*(M_N(\bbA^\hbar[\hbar^{-1}])[\eta ])$ is very homotopically constant if we replace $\iota_{\underline{g}}$ and $L_{\underline{g}}$ by $\iota_{-\underline{g}}$ and $L_{-\underline{g}}$. )

As in the case of very homotopically constant modules we get a morphism of complexes 
$$\overline{C}^\lambda_*(M_N(\bbA^\hbar [\hbar^{-1}])[\eta ])_{\fh
  +\underline{\fh}} \rightarrow C^*( \tilde{\fh}[\epsilon ], \tilde{\fh};\overline{C}^\lambda_*(M_N(\bbA^\hbar [\hbar^{-1}])[\eta ])$$
given by 
$$l \rightarrow ((h_1\epsilon, \ldots , h_p\epsilon ) \rightarrow (-1)^p\iota_{\underline{h}_1} \cdots \iota_{\underline{h}_p}l).$$

Note that $\eta^{(k)} \in \overline{C}^\lambda_*(M_N(\bbA^\hbar [\hbar^{-1}])[\eta ])_{\fh
  +\underline{\fh}}$ and we therefore get classes
$$\eta^{[k]}\hbox{ in } C^*(\tilde{\fh}[\epsilon
],\tilde{\fh}; \overline{C}^\lambda_*(M_N(\bbA^\hbar
[\hbar^{-1}])[\eta]))$$
given by
$$(h_1\epsilon, \ldots , h_p\epsilon )\rightarrow
(-1)^p\iota_{h_1\eta}\cdots \iota_{h_p\eta}\eta^{(k)}.$$

\begin{lemma}  In $H^*(\tilde{\fh}[\epsilon ],\tilde{\fh};
\overline{C}^\lambda_*(M_N(\bbA^\hbar[\hbar^{-1}])[\eta]))$ one has the formula
$$U=\sum_{m \geq 0}(\hat{A}\cdot e^{-\theta} \cdot ch^{-1})^{-1}_{2m}\cdot
\eta^{[m]}.$$
\end{lemma}

\textit{Proof.} It is well known from \cite{BNT} that there is a
splitting principle, i.e. the inclusion morphism
\begin{eqnarray*}
& H^*(\tilde{\fh} [\epsilon ],\tilde{\fh};
\overline{C}^\lambda_*(M_N(\bbA^\hbar[\hbar^{-1}])[\eta]))
\rightarrow \\
&
H^*((\tilde{\fd}_n+\fsu(N)) [\epsilon ],(\tilde{\fd}_n+\fsu(N));
\overline{C}^\lambda_*(M_N(\bbA^\hbar[\hbar^{-1}])[\eta]))
\end{eqnarray*}
where  $\tilde{\fd}_n=\fd_n+\hbar^{-1}\bbC[[\hbar ]]$ and $\fd_n$  is
the set of $n\times n$ diagonal matrices, is injective.

Therefore, we only have to identify the two
classes in  $$H^*((\tilde{\fd}_n+\fsu(N)) [\epsilon ],(\tilde{\fd}_n+\fsu(N));
\overline{C}^\lambda_*(M_N(\bbA^\hbar[\hbar^{-1}])[\eta])).$$

We next note that we can factor the classes at hand in the following
way: Write
$$M_N(\bbA^\hbar[\hbar^{-1}])=\bbA^\hbar_1[\hbar^{-1}]\otimes \ldots
\otimes \bbA^\hbar_1[\hbar^{-1}]\otimes M_N(\bbA^\hbar_1) [\hbar^{-1}],$$
where $\bbA^{\hbar}_1$ is the formal Weyl algebra in one variable. We can write $U=U_1\times \ldots \times U_1\times U_1'$, where
$U_1$ is the extension of the fundamental class in
$C^*(\fd_1[\epsilon],\fd_1;\overline{C}^\lambda_*(\bbA^\hbar_1[\hbar^{-1},\eta]))$
and $U_1'$ is the extension of the fundamental class in
$$
C^*((\tilde{\fd}_1+\fsu(N))[\epsilon],\tilde{\fd}_1+\fsu(N);\overline{C}^\lambda_*(M_N(\bbA^\hbar_1[\hbar^{-1}])[\eta])).$$
We thus need to identify $U_1$ and $U_1'$.

In the first case we can represent $U_1$ by
\begin{equation} \label{fund}
U_1=\sum_{m=1}^\infty\frac{1}{m}((i\hbar)^{-1}\hat{\xi}\otimes
\hat{x})^{\otimes m}c_1^{m-1}
\end{equation}
where $c_1$ is the first Chern class.

Recall that the definition of $\bbA^\hbar_1[\hbar^{-1}]$ is $\bbC
[[\hat{x}, \hat{\xi} ]][[\hbar, \hbar^{-1}]$ with a product
$*$. Given an element $f$ in $\bbA^{\hbar}_1[\hbar^{-1}]$, we can regard
$f$ as a function in the variables $\hat{x}, \hat{\xi}$ with values in
$\bbC[[\hbar, \hbar^{-1} ]$. Hence, given $f$ in $\bbA^\hbar_1
[\hbar^{-1}]$, we can define $l( f)=f(0,0)$. With this $l$ we get,
according to section \ref{cyclic}, a
quasi-isomorphism of complexes
$$Br : \overline{C}^\lambda_*(\bbA^\hbar [\hbar^{-1}, \eta ])
\rightarrow
C^\lambda_*(\bbC [[\hbar, \hbar^{-1}]).$$

One checks that this gives a morphism of complexes
\begin{eqnarray*}  &C^*(\fd_1[\epsilon ],\fd_1;\overline{C}_*^\lambda
(\bbA^\hbar_1[\hbar^{-1},\eta]))
 \stackrel{Br}{\longrightarrow}
 C^*(\fd_1[\epsilon ], \fd_1 ; C^\lambda_*(\bbC[[\hbar,\hbar^{-1}]
))
\end{eqnarray*}
and therefore a quasi isomorphism of complexes.

A computation now shows that
$$Br(U_1)= \sum_{m=0}^\infty  1^{(m)}\hat{A}^{-1}_{2m},$$
where $1^{(m)}=m!(m+1)!1^{\otimes (2m+1)}$ and $\hat{A}$ is as in
example 3 in section \ref{liesekt}. On the other hand we have
$Br(\eta^{[m]})=1^{(m)}$, where $\eta^{[m]}$ is the class
$$\epsilon d_1, \ldots , \epsilon d_k \rightarrow (-1)^k\iota_{d_1\eta }\cdots
\iota_{d_k \eta}\eta^{(m)}, \quad d_i \in \fd_1.$$
Since $Br$ is a quasi-isomorphism, we have
$$U_1 = \sum_{m=0}^\infty \eta^{[m]}\hat{A}^{-1}_{2m}$$
in $H^*(\fd_1[\epsilon ],\fd_1 ;
\hat{C}^\lambda_*(\bbA_1^\hbar[\hbar^{-1}, \eta ]))$.

We note that we have a morphism of complexes
\begin{eqnarray*}
&C^*((\tilde{\fd}_1+\fsu(N))[\epsilon],\tilde{\fd}_1+\fsu(N);\overline{C}^\lambda_*(M_N(\bbA^\hbar[\hbar^{-1}])))\\
& Tr \downarrow\\
&
C^*((\tilde{\fd}_1+\fsu(N))[\epsilon],\tilde{\fd}_1+\fsu (N);\overline{C}^\lambda_*(\bbA^\hbar
[\hbar^{-1}]))
\end{eqnarray*}
where in the bottom row $\fsu (N)$ acts trivially and $Tr$ denotes the
morphism of complexes
$$\overline{C}^\lambda_*(M_N(\bbA^\hbar_1)[\hbar^{-1}]) \rightarrow
\overline{C}^\lambda_*(\bbA^\hbar_1[\hbar^{-1}])$$
induced by the normalized trace $Tr$.

In the top row we
have the class $U_1'$, that in homology is the unique extension of the
fundamental class.
In the bottom row we have a class $U'$, given by the same formula as in
\ref{fund}. Also this is, in homology, a
unique extension of the fundamental class. Therefore $Tr (U)= U'$ in
$$
H^*((\tilde{\fd}_1+\fsu(N)[\epsilon],\tilde{\fd}_1+\fsu(N);\overline{C}^\lambda_*(\bbA^\hbar[\hbar^{-1}])).$$
We further note that there are morphisms of complexes
\begin{eqnarray*}
&C^*((\tilde{\fd}_1+\fsu(N))[\epsilon],\tilde{\fd}_1+\fsu(N);\overline{C}^\lambda_*(M_N(\bbA^\hbar[\hbar^{-1}])[\eta]))\\
& Tr \downarrow \\
&
C^*((\tilde{\fd}_1+\fsu(N))[\epsilon],\tilde{\fd}_1+\fsu(N);\overline{C}^\lambda_*(\bbA^\hbar
[\hbar^{-1}][\eta]))
\\
& Br \downarrow \\
& C^*((\tilde{\fd}_1+\fsu(N))[\epsilon],\tilde{\fd}_1+\fsu(N);C^\lambda_*(\bbC[[\hbar, \hbar^{-1}]))
\end{eqnarray*}
where $Br$ is the Brodzki map as before. We thus have
$$Br( Tr (U_1'))=\sum_{m \leq 0} \hat{A}^{-1}_{2m}1^{(m)}.$$
Note that
$$Br(Tr(\eta^{[m]}))= \sum_{l=0}^\infty 1^{(m+l)}(e^\theta ch)_{2l}$$
where $\theta$ is given in example 1 in section \ref{liesekt} and
$ch$ is given in example 2 in section \ref{liesekt}.

Since $Br\circ Tr$ is a quasi-isomophism, we get
$$U_1'=\sum_{m=0}^\infty (\hat{A}\cdot e^{-\theta} \cdot (ch)^{-1})^{-1}_{2m}\eta^{[m]}$$
in
$H^*((\tilde{\fd}_1+\fsu(N))[\epsilon],\tilde{\fd}_1+\fsu(N);\overline{C}^\lambda_*(M_N(\bbA^\hbar_1[\hbar^{-1}])[\eta]))$.

The lemma now follows, since $\hat{A}$ is multiplicative.
\\

It is well known, see \cite{BNT}, that the restriction homomorphism
\begin{eqnarray*} & C^*(\tilde{\fg}[\epsilon
],\tilde{\fh}; \overline{C}^\lambda_*(M_N(\bbA^\hbar
[\hbar^{-1}])[\eta])) \rightarrow
 C^*(\tilde{\fh}[\epsilon
],\tilde{\fh}; \overline{C}^\lambda_*(M_N(\bbA^\hbar
[\hbar^{-1}])[\eta]))
\end{eqnarray*}
is a quasi-isomorphism. But in $ C^*(\tilde{\fg}[\epsilon
],\tilde{\fh}; \overline{C}^\lambda_*(M_N(\bbA^\hbar
[\hbar^{-1}])[\eta]))$ there are two classes that maps to $\eta^{[m]}$  under the
restriction homomorphism, namely $CW(\eta^{[m]})$
and the class
\begin{eqnarray}\label{CW}
(g_1\epsilon, \ldots , g_p\epsilon ) \rightarrow
(-1)^p\iota_{g_1\eta}\cdots \iota_{g_p\eta}
\eta^{(m)}.
\end{eqnarray}
Therefore
$CW(\eta^{[m]})$ and the class \ref{CW} are equivalent in $$ C^*(\tilde{\fg}[\epsilon
],\tilde{\fh}; \overline{C}^\lambda_*(M_N(\bbA^\hbar
[\hbar^{-1}])[\eta])).$$
Considering the restriction homomorphism
$$ C^*(\tilde{\fg}[\epsilon
],\tilde{\fh}; \overline{C}^\lambda_*(M_N(\bbA^\hbar
[\hbar^{-1}])[\eta])) \rightarrow
 C^*(\tilde{\fg},\tilde{\fh}; \overline{C}^\lambda_*(M_N(\bbA^\hbar
[\hbar^{-1}])[\eta]))$$
we get that in the righthand side, $CW(\eta^{[m]})$,
is equivalent to $\eta^{(m)}$. For $\fg$-modules $\bbL^*$ we
have that $C^*(\tilde{\fg},\tilde{\fh} ; \bbL^*)$ is quasi-isomorphic to
$C^*(\fg,\fh;\bbL^*)$. Combining these observations we have

\begin{thm}  \label{teks}
Let $U$ be an extension of $U_0$ to a class in  $$C^*(\fg
,\fh; \overline{C}^\lambda_*(M_N(\bbA^\hbar
[\hbar^{-1}])[\eta])).$$ Then in
 $H^*(\fg
,\fh; \overline{C}^\lambda_*(M_N(\bbA^\hbar
[\hbar^{-1}])[\eta]))$ we have the following equality

$$U=\sum_{m \geq 0}(\hat{A}\cdot e^\theta \cdot (ch)^{-1})^{-1}_{2m}\cdot
\eta^{(m)}.$$
\end{thm}

\section{The Gelfand-Fuks Construction}
We now consider a $\fg$-module $\bbL^*$, where $\fg$ is as in section \ref{defqe}. Given a deformation
quantization $A^\hbar_E$ we can consider the bundle
$$\tilde{M}_{A^\hbar_E}\times_G\bbL^*$$
and also consider the differential forms with values in this
bundle. We will denote this by $\Omega^*(M,\bbL^*)$. Furthermore we get
a flat connection $\nabla$ induced from the connection on
$\tilde{M}_{A^\hbar_E}$. Using this connection and the differential on
$\bbL^*$ we  get a complex
$\Omega(M,\bbL^*)$. The Gelfand-Fuks construction gives a morphism of
complexes
$$GF:C^*(\fg,\fh;\bbL^*) \rightarrow \Omega^*(M,\bbL^*)$$
defined in the following way:

Choose a $U(N)\times SU(N)$ trivialization of
$\tilde{M}_{A^\hbar_E}\times_G\bbL^*$. In a given
trivialization  write $\nabla=d+A$, where $A$
is the connection one form. Given vector fields $X_1, \ldots , X_p$
 and $l$ in $C^*(\fg,\fh ; \bbL^*)$  define
$$GF(l)(X_1, \ldots , X_p )=l(X_1, \ldots X_p ).$$

If we look at the examples of classes in $C^*(\fg, \fh; \bbL^*)$ 
constructed in section \ref{liesekt}, we see that for $\theta$ in
$C^*(\fg , \fh ; \frac{1}{\hbar}\bbC [[ \hbar ]])$ we get that
$GF(\theta )$ is the characteristic class of the deformation
quantization.
In the case of $\hat{A}$ we see that $GF (\hat{A})$ is the $\hat{A}$
class of $TM$, and the case of $ch$ this is just the Chern character of
$End(E)$.

The main example we are going to look at is the case where $\bbL^*$ is
$\overline{C}^\lambda_*(M_N(\bbA^\hbar)[\hbar^{-1}])$. We first note that
\begin{lemma} \label{acyk}
Let $A_E^\hbar$ be a deformation quantization over $\bbR^{2n}$.
The complex
$(\Omega^*(\bbR^{2n},\overline{C}^\lambda_*(M_N(\bbA^\hbar))),\nabla)$ is
acyclic and the cohomology consists of jets on the diagonal of elements in
$\overline{C}^\lambda_*(M_N(\cW_n))$.
\end{lemma}

\textit{Proof.} We can assume that $\nabla $ is of the form
$d-\sum_i(\partial_{\hat{x}_i}\otimes dx_i
+\partial_{\hat{\xi}_i}\otimes d\xi_i)$. We have a short exact
sequence of complexes
\begin{eqnarray*}
0 \rightarrow
(\Omega^*(\bbR^{2kn}/\Delta,M_N(\bbA^\hbar)^{\otimes k}),\nabla)
\rightarrow
(\Omega^*(\bbR^{2kn},M_N(\bbA^\hbar)^{\otimes k}),\nabla)
\stackrel{\varphi^*}{\rightarrow}
\\
(\Omega^*(\bbR^{2k},M_N(\bbA^\hbar)^{\otimes k}), \nabla) \rightarrow 0
\end{eqnarray*}
where $\varphi :\bbR^{2n} \rightarrow \bbR^{2kn}$ is the map onto the
diagonal $\delta$. From the associated long exact sequence we  see that
$(\Omega^*(\bbR^{2k},M_N(\bbA^\hbar)^{\otimes k}), \nabla)$ is acyclic
and that the cohomology consists of jets on the diagonal of elements of
$M_N(\cW_n)^{\otimes k}$. By considering the following short exact
sequence of complexes
\begin{eqnarray*}
& 0 \rightarrow
(\Omega^*(\bbR^{2k},1 \otimes M_N(\bbA^\hbar)\otimes \ldots \otimes
M_N(\bbA^\hbar)+\ldots \\& +M_N(\bbA^\hbar)\otimes \ldots \otimes
M_N(\bbA^\hbar)\otimes 1), \nabla) \\
&  \rightarrow  (\Omega^*(\bbR^{2k},M_N(\bbA^\hbar)^{\otimes k}),
 \nabla) \rightarrow
 (\Omega^*(\bbR^{2k},\overline{M_N(\bbA^\hbar)}^{\otimes k}), \nabla)
 \rightarrow  0
\end{eqnarray*}
we  see that $(\Omega^*(\bbR^{2k},\overline{M_N(\bbA^\hbar)}^{\otimes
  k}), \nabla)$ is acyclic and that the cohomology consists of jets on the
diagonal of elements in $\overline{M_N(\cW_n)}^{\otimes k}$.

Let  $a$ be an element  in
$\Omega^*(M,\overline{M_N(\bbA^\hbar)}^{\otimes k}/Im(1-\tau ))$ with
$\nabla (a)=0$. We can then lift $a$ to an element $\tilde{a} \in \Omega^*(\bbR^{2k},\overline{M_N(\bbA^\hbar)}^{\otimes
  k})$, where $\nabla (\tilde{a}) \in \Omega^*(\bbR^{2k},Im(1-\tau
))$. However, $b= \frac{1}{k}\sum_{i=0}^{k-1} \tau^i(\tilde{a})$ is also a lift of $a$
and $\nabla (b)=\sum_{i=0}^{k-1} \tau^i \nabla ( \tilde{a})=0$. The
lemma follows from this.
\\

\section{Traces on Deformation Quantizations and Index Theory}

We consider a deformation quantization $A^\hbar_E$ of an endomorphism bundle
$End(E)$ over a symplectic manifold $M$ of dimension $n$. Let $A^\hbar_{E,c}$ be the
algebra of elements in $A^\hbar_E$ with compact support. This algebra has a canonical
$\bbC [[\hbar, \hbar^{-1}]$ valued trace defined in the following way:

Let $(V_i,\Phi_i)$ be a cover of $M$ and
$\Phi_i^{-1}:M_N(\cW_n)\rightarrow A^\hbar_E$ local
isomorphisms over $V_i$. Let $\rho_{V_i}$ be a partition of unity with
respect to the covering. For an element $a \in A^\hbar_{E,c}$ we define
$$ Tr(a)=\sum_i \int \frac{1}{n!(i\hbar)^n}tr(\Phi_i(\rho_i *
a))\omega_{st}^n,$$ where $\omega_{st}$ is the standard symplectic form
on $\bbR^{2n}$ and $tr$ denotes the normalized trace on $M_N$. That $Tr$ is independent of the choices made, and that it is a trace,
hinges on the following two propositions.
\begin{proposition}
Let $E$ be the trivial line bundle. Then $Tr$ is a trace and
independent of the choices made.
\end{proposition}

\textit{Proof.} See \cite{Fe1}.
\\

\begin{proposition} If $\cW_n$ is the Weyl algebra over some
  contractible open subset $U$ of  $\bbR^{2n}$, then every automorphism over the
  identity map of $M_N(\cW_n)$ is inner.
\end{proposition}

\textit{Proof.} More or less the same as lemma \ref{inner}, see also
lemma \ref{cont}.
\\

We consider the trace as a functional on $CC^{per}_*(A^\hbar_{E,c})$ and
we want to compute the trace at the level of homology. To this end we
consider the following situation.

Given an element $b$ in
$\overline{C}^\lambda_*(A^\hbar_E[\hbar^{-1},\eta])$, we define
$\chi_{Tr}(b)$ in $CC^*_{per}(A^\hbar_{E,c})$ in the following way
$$\chi_{Tr}(b)(a)=Tr( b \cdot a ),$$
where $ \cdot $ means the action of
$\overline{C}^\lambda_*(A^\hbar_E[\hbar^{-1},\eta])$ on
$CC^{per}_*(A_{E,c}^\hbar[\hbar^{-1}])$, see Theorem \ref{opera}. We will now extend this to
elements in the \v{C}ech complex
$\check{C}^*(M,\overline{C}^\lambda_*(A^\hbar_E[\hbar^{-1},\eta])$
with values in the presheaf $V \rightarrow
\overline{C}^\lambda_*(A^\hbar_{E|V}[\hbar^{-1},\eta])$. This is done
in the following proposition.
\begin{proposition}
Take $\{ b_{V_0 \ldots V_p}\}$ in
$\check{C}^*(M,\overline{C}^\lambda_*(A^\hbar_E[\hbar^{-1},\eta]))$ and
$a$ in $CC_*^{per}(A^\hbar_{E,c}[\hbar^{-1}])$. Define
$$\chi_{Tr}(\{ b_{V_0 \ldots V_p}\} )(a)=\sum_{V_0 , \ldots ,
  V_p}\chi_{Tr}(b_{V_0, \ldots , v_p})(I_{\rho _{V_0}}[B+b,I_{\rho_{ V_1}}]\ldots [B+b,I_{\rho_{ V_p}}]a).
$$
This gives a morphism of complexes
$$\check{C}^*(M,\overline{C}^\lambda_*(A^\hbar_E[\hbar^{-1},\eta]))
\otimes  CC_*^{per}(A^\hbar_{E,c}[\hbar^{-1}]) \rightarrow \bbC
[[\hbar, \hbar^{-1}].$$

\end{proposition}

\textit{Proof.} See \cite{NT1}.

\subsection{The Fundamental Class in the \v{C}ech Complex}
Recall that we have the canonical coordinates $x_1, \ldots , x_n ,
\xi_1, \ldots , \xi_n $ on $\bbR^{2n}$. We will also use this notation
for the associated
coordinate functions and consider these coordinate
functions as elements in $\cW_n$.

We can consider the fundamental  class $U_0$ in
$\overline{C}^\lambda_*(M_N(\cW_n)[\hbar^{-1}])$ given by
$$U_0=\frac{1}{2n(i\hbar)^n }\sum_{\sigma \in
  S_{2n}}(v_{\sigma_{1}}\otimes \ldots \otimes v_{\sigma_{2n}})$$
where $(v_1, \ldots , v_{2n})=(x_1, \ldots , x_n, \xi_1, \ldots ,
\xi_n )$.

By the same argument as in the section on the fundamental class in  Lie-algebra cohomology,
this class extends uniquely in cohomology to a class $U$ in
$\check{C}^*(M,\overline{C}^\lambda_*(A^\hbar_E[\hbar^{-1}]))$. In order to connect
this class to the fundamental class defined in Lie algebra cohomology
we introduce the complex
$$\check{C}^*(M,\Omega^*(M,\overline{C}^\lambda_*(M_N(\bbA^\hbar[\hbar^{-1}]))).$$
According to Lemma \ref{acyk} this is quasi-isomorphic to the complex
$\check{C}^*(M,\overline{C}^\lambda_*(A^\hbar_E[\hbar^{-1}]))$.
Furthermore we have the Gelfand-Fuks morphism
$$GF:C^*(\fg,\fh;\overline{C}^\lambda_*(M_N(\bbA^\hbar[\hbar^{-1}])) \rightarrow \check{C}^*(M,\Omega^*(M,\overline{C}^\lambda_*(M_N(\bbA^\hbar[\hbar^{-1}]))).$$
Because of uniqueness we get that 
$GF(U)=U$ in cohomology. By theorem \ref{teks} and the splitting principle we therefore get
\begin{thm} \label{fundamen}
In the complex
$\check{C}^*(M,\overline{C}^\lambda_*(A^\hbar_E[\hbar^{-1},\eta]))$
the two classes 
$$U \quad \hbox{ and } \quad \sum_{m \geq 0}(\hat{A}\cdot e^\theta \cdot ch)^{-1}_{2m}\cdot
\eta^{(m)}$$
are equivalent.
\end{thm}
With this we are now in position to prove
\begin{thm} \label{hoved}
 $\chi_{Tr}(U)(a_0 \otimes \ldots \otimes a_k )$ has no
singularities in $\hbar$ and
$$\chi_{Tr}(U)(a_0 \otimes \ldots \otimes a_k)=(-1)^n\int
ch^{-1}(End(E))ch(\nabla)(\tilde{a}_0\otimes \ldots \otimes  \tilde{a}_k) \hbox{ mod
 }\hbar . $$
Here $\tilde{a}_i$ is $a_i$ mod $\hbar$ and $ch(\nabla )$ is the
J.L.O. cocycle associated to $\nabla$ (see also \cite{GO}), i.e.
\begin{eqnarray*} && ch(\nabla)(\tilde{a}_0 \otimes \ldots \otimes
  \tilde{a}_k)= \\
&& \int_{\Delta_k}tr (\tilde{a}_0e^{-t_0\nabla^2}\nabla
(\tilde{a}_i)e^{-t_1\nabla^2}\cdots \nabla
(\tilde{a}_k)e^{-t_k\nabla^2} )dt_0 \cdots dt_{k-1}
\end{eqnarray*}
where $tr$ is the normalized trace on $End(E)$.

 According to Theorem \ref{fundamen} we therefore have that
$$Tr(a_0 \otimes \ldots \otimes a_k \cdot e^{-\theta} )=(-1)^n\int
\hat{A} \cdot ch(\nabla )(\tilde{a}_0\otimes \ldots \otimes \tilde{a}_k)
\hbox{ mod } \hbar$$
\end{thm}

\textit{Proof.} Because of Morita equivalence it is enough to look at
the case where $\tilde{a}_i$ is scalar for all $i$. We  have that
$$\chi_{Tr}(U)(a_0 \otimes \ldots \otimes a_k)=
\sum_{V_0}\chi_{tr}(U_0)(I_{\rho V_0}(a_0\otimes \ldots \otimes a_k))+\ldots$$
and it is not difficult to see that $\ldots$ is zero modulo
$\hbar$. The explicit formula for $\chi_{Tr}(U_0)$ gives
$$\sum_{V_0}\chi_{Tr}(U_0)(I_{\rho V_0}(a_0 \otimes \ldots \otimes a_{2n}))=(-1)^n\frac{1}{2n!}\int
\tilde{a}_0d\tilde{a}_1 \cdots d \tilde{a}_{2n} \hbox{ mod }\hbar.$$
The result follows from this.\\

By adopting the arguments in \cite{NT5} one sees that 
$$\frac{d}{d\hbar }Tr (a \cdot e^{-\theta } )=0$$
when $a$ is a cycle in $CC^{per, \bbC}_*(A°^\hbar_{E,c})$. Here $CC^{per,\bbC}_*(A^\hbar_{E,c})$ means cyclic periodic homology of $A^\hbar_{E,c}$ as a $\bbC$-algebra. 
Together with theorem \ref{hoved} we get
\begin{thm}
The identity 
$$Tr(a)=(-1)^n\int \hat{A}\cdot e^\theta \cdot ch(\nabla)(\tilde{a})$$
holds when $a$ is a cycle in $CC^{per, \bbC}_*(A°^\hbar_{E,c})$ 
\end{thm}

\bibliographystyle{alpha}
\bibliography{ref}
\end{document}